\documentclass{amsart}

\usepackage{latexsym}

\usepackage{amsfonts,amssymb,euscript,epsfig}

\usepackage[all]{xy}

\newtheorem{theo}{Theorem}[section]
\newtheorem{defi}[theo]{Definition}
\newtheorem{propdefi}[theo]{Proposition-Definition}

\newtheorem{lem}[theo]{Lemma}
\newtheorem{prop}[theo]{Proposition}
\newtheorem{rem}[theo]{Remark}
\newtheorem{coro}[theo]{Corollary}

\newcommand{\hgot}{\ensuremath{\mathfrak{h}}}
\newcommand{\kgot}{\ensuremath{\mathfrak{k}}}

\newcommand{\tgot}{\ensuremath{\mathfrak{t}}}
\newcommand{\ugot}{\ensuremath{\mathfrak{u}}}

\newcommand{\Acal}{\ensuremath{\mathcal{A}}}

\newcommand{\Ccal}{\ensuremath{\mathcal{C}}}

\newcommand{\Fcal}{\ensuremath{\mathcal{F}}}

\newcommand{\Kcal}{\ensuremath{\mathcal{K}}}
\newcommand{\Lcal}{\ensuremath{\mathcal{L}}}

\newcommand{\Ocal}{\ensuremath{\mathcal{O}}}
\newcommand{\Qcal}{\ensuremath{\mathcal{Q}}}

\newcommand{\Xcal}{\ensuremath{\mathcal{X}}}
\newcommand{\Ycal}{\ensuremath{\mathcal{Y}}}
\newcommand{\Zcal}{\ensuremath{\mathcal{Z}}}
\newcommand{\Ucal}{\ensuremath{\mathcal{U}}}
\newcommand{\Vcal}{\ensuremath{\mathcal{V}}}

\newcommand{\Pbb}{\ensuremath{\mathbb{P}}}
\newcommand{\Z}{\ensuremath{\mathbb{Z}}}
\newcommand{\C}{\ensuremath{\mathbb{C}}}

\newcommand{\R}{\ensuremath{\mathbb{R}}}
\newcommand{\N}{\ensuremath{\mathbb{N}}}

\newcommand{\croc}{\ensuremath{\hookrightarrow}}

\newcommand{\indice}{\ensuremath{\hbox{\rm Index}}}
\newcommand{\image}{\ensuremath{\hbox{\rm Image}}}

\newcommand{\Cr}{\ensuremath{\hbox{\rm Cr}}}

\newcommand{\Tr}{\ensuremath{\hbox{\rm Tr}}}

\newcommand{\Char}{\ensuremath{\hbox{\rm Char}}}
\newcommand{\End}{\ensuremath{\hbox{\rm End}}}

\newcommand{\qfor}{\ensuremath{\mathcal{Q}^{-\infty}}}
\newcommand{\Rfor}{\ensuremath{R^{-\infty}}}
\newcommand{\GL}{\ensuremath{GL}}

\def \K {{\rm \bf K}}
\def \T {{\rm T}}
\def \h {{\rm h}}
\def \what {\widehat}

\def \U {{\rm U}}
\def \ad {{\rm ad}}

\def \clif {{\bf c}}

\title{Formal Geometric Quantization}

\author{Paul-Emile  PARADAN}

\address{Institut de Math\'ematiques et Mod\'elisation de Montpellier (I3M), Universit\'e Montpellier 2} 

\email{Paul-Emile.Paradan@math.univ-montp2.fr}

\date{8 February 2007}

\begin{document}


\begin{abstract}
Let $K$ be a compact Lie group acting in an Hamiltonian way on a
symplectic manifold $(M,\Omega)$ which is prequantized by a
Kostant-Souriau line bundle. We suppose here that the moment map
$\Phi$ is \emph{proper} so that the reduced space
$M_{\mu}:=\Phi^{-1}(K\cdot\mu)/K$ are compact for all $\mu$.
Following Weitsman \cite{Weitsman}, we define the ``formal geometric
quantization'' of $M$ as
$$
\qfor_K(M):=\sum_{\mu\in \what{K}} \Qcal(M_{\mu}) V_\mu^K.
$$
The aim of this article is to study the functorial properties of
the assigment $M\longrightarrow \qfor_K(M)$.
\end{abstract}


\maketitle


{\small
\tableofcontents}

\section{Introduction and statement of the results}\label{sec:intro}

Let $M$ be an Hamiltonian $K$-manifold with symplectic form $\Omega$
and moment map $\Phi: M\to \kgot^{*}$. We assume the existence of a
$K$-equivariant line bundle $L$ on $M$ with connection of curvature
equal to $-i\, \omega$. In other words $M$ is pre-quantizable in the
sense of \cite{Kostant70} and we call $L$ the Kostant-Souriau line
bundle.

In the process of quantization one tries to associate a unitary
representation of $K$ to  the data $(M,\Omega,\Phi, L)$. When $M$ is
\textbf{compact} one associates to these data a virtual
representation $\Qcal_K(M)\in R(K)$ of $K$ defined as the
equivariant index of a Dolbeault-Dirac operator : $\Qcal_K(M)$ is
the geometric quantization of $M$.

This quantization process satisfies the following functorial
properties :

$\textbf{[P1]}$  When $N$ and $M$ are respectively pre-quantized
compact Hamiltonian $K_1$ and $K_2$ manifolds, the product $M\times
N$ is a pre-quantized compact Hamiltonian $K_1\times K_2$-manifold,
and we have
\begin{equation}\label{eq:functor-produit}
\Qcal_{K_1\times K_2}(M\times N)=\Qcal_{K_1}(M)\otimes
\Qcal_{K_2}(N) \quad \mathrm{in} \quad R(K_1\times K_2)\simeq
R(K_1)\otimes R(K_2).
\end{equation}

\medskip

$\textbf{[P2]}$ If $H\subset K$ is a connected Lie subgroup, then
the restriction of $\Qcal_K(M)$ to $H$ is equal to $\Qcal_H(M)$.

\medskip

Note that $\textbf{[P1]}$ and $\textbf{[P2]}$ gives the following
functorial property :

\medskip

$\textbf{[P3]}$  When $N$ and $M$ are pre-quantized compact
Hamiltonian $K$-manifold, the product $M\times N$ is a pre-quantized
compact Hamiltonian $K$-manifold, and we have $\Qcal_{K}(M\times
N)=\Qcal_{K}(M)\cdot\Qcal_{K}(N)$, where $\cdot$ denotes the product
in $R(K)$.

\medskip

One other fundamental property is the behaviour of the
$K$-multiplicities of $\Qcal_K(M)$ that is known as ``quantization
commutes with reduction''.

Let $T$ be a maximal torus of $K$, $C_K\subset \tgot^*$ be a Weyl
Chamber, $\wedge^*\subset\tgot^*$ be the weight lattice, and
$\what{K}=\wedge^*\cap C_K$ be the set of dominant weights. The ring
of character $R(K)$ as a $\Z$-basis $V_\mu^K,\mu\in \what{K}$ :
$V_\mu^K$ is the irreducible $K$-representation with highest weight
$\mu$.

For any $\mu\in \what{K}$ which is a regular value of $\Phi$, the
reduced space (or symplectic quotient)
$M_{\mu}:=\Phi^{-1}(K\cdot\mu)/K$ is an orbifold equipped with a
symplectic structure $\Omega_{\mu}$. Moreover
$L_{\mu}:=(L\vert_{\Phi^{-1}(\mu)}\otimes \C_{-\mu})/K_{\mu}$ is a
Kostant line orbibundle over $(M_{\mu},\Omega_{\mu})$. The
definition of the index of the Dolbeault-Dirac operator carries over
to the orbifold case, hence $\Qcal(M_{\mu})\in \Z$ is defined. In
\cite{Meinrenken-Sjamaar}, this is extended further to the case of
singular symplectic quotients, using partial (or shift)
de-singularization. So the integer $\Qcal(M_{\mu})\in \Z$ is well
defined for every $\mu\in \what{K}$ : in particular
$\Qcal(M_{\mu})=0$ if $\mu\notin\Phi(M)$.

The following Theorem was conjectured by Guillemin-Sternberg
\cite{Guillemin-Sternberg82} and is known as ``quantization commutes
with reduction''
\cite{Meinrenken98,Meinrenken-Sjamaar,Tian-Zhang98,pep-RR}. For
complete references on the subject the reader should consult
\cite{Sjamaar96,Vergne-Bourbaki}.

\begin{theo} (Meinrenken, Meinrenken-Sjamaar). \label{theo:Q-R}
We have the following equality in $R(K)$
$$
\Qcal_K(M)=\sum_{\mu\in \what{K}}\Qcal(M_{\mu})\, V_{\mu}^{K}\ .
$$
\end{theo}

\medskip

Suppose now that $M$ is \textbf{non-compact} but the moment map
$\Phi: M\to\kgot^*$ is assumed to be \textbf{proper} (we will said
simply "$M$ is proper"). In this situation the geometric
quantization of $M$ as an index of an elliptic operator is not well
defined. Nevertheless the integers $\Qcal(M_{\mu}),\mu\in \what{K}$
are well defined since the symplectic quotient $M_\mu$ are
\textbf{compact}.

Following Weitsman \cite{Weitsman}, we introduce the following
\begin{defi}\label{def:formal-quant}
The formal quantization of $(M,\Omega,\Phi)$ is the element of
$\Rfor(K):=\hom_\Z(R(K),\Z)$ defined by
$$
\qfor_K(M)=\sum_{\mu\in \what{K}}\Qcal(M_{\mu})\, V_{\mu}^{K}\ .
$$
\end{defi}

\medskip

A representation $E$ of $K$ is  admissible if it as finite
$K$-multiplicities : \break $\dim(\hom_K(V_\mu^K,E))<\infty$ for
every $\mu\in\what{K}$. Here $\Rfor(K)$ is the Grothendieck group
associated to the $K$-\emph{admissible} representations. We have a
canonical inclusion $i:R(K)\croc \Rfor(K)$ : to $V\in R(K)$ we
associate the map $i(V):R(K)\to \Z$ defined by $W\mapsto
\dim(\hom_K(V,W))$. In order to simplify the notation, $i(V)$ will
be written $V$. Moreover the tensor product induces a structure of
$R(K)$-module on $\Rfor(K)$ since $E\otimes V$ is an admissible
representation when $V$ and $E$ are respectively finite dimensional
and admissible representation of $K$.

\medskip

It is an easy matter to see that $\textbf{[P1]}$ holds for the
formal quantization process $\qfor$. Let $N$ and $M$ be respectively
pre-quantized \emph{proper} Hamiltonian $K_1$ and $K_2$ manifolds :
the product $M\times N$ is then a pre-quantized \emph{proper}
Hamiltonian $K_1\times K_2$-manifold. For the reduced spaces we have
$(M\times N)_{(\mu_1,\mu_2)}\simeq M_{\mu_1} \times N_{\mu_2}$, for
all $\mu_1\in\what{K}_1$, $\mu_2\in\what{K}_2$. It follows then that
\begin{equation}\label{eq:functor-produit-formal}
\qfor_{K_1\times K_2}(M\times N)=\qfor_{K_1}(M)\what{\otimes}\,
\qfor_{K_2}(N)
\end{equation}
in $\Rfor(K_1\times K_2)\simeq \Rfor(K_1)\what{\otimes} \Rfor(K_2)$.

\medskip

The purpose of this article is to show that the functorial property
$\textbf{[P2]}$ still holds for the formal quantization process
$\qfor$.

\begin{theo}\label{theo:intro}
    Let $M$ be a pre-quantized
Hamiltonian $K$-manifold which is \emph{proper}. Let $H\subset K$ be
a connected Lie subgroup such that $M$ is still \emph{proper} as an
Hamiltonian $H$-manifold. Then $\qfor_K(M)$ is $H$-admissible and we
have the following equality in $\Rfor(H)$ :
\begin{equation}\label{eq:functor-restriction-formal}
\qfor_K(M)|_H=\qfor_H(M).
\end{equation}
\end{theo}

\medskip

For $\mu\in \what{K}$ and $\nu\in \what{H}$ we denote
$N^\mu_{\nu}=\dim (\hom_H(V_\nu^H,V_\mu^K|_H))$ the multiplicity of
$V_\nu^H$ in the restriction $V_\mu^K|_H$. In the situation of
Theorem \ref{theo:intro}, the moment maps relative to the $K$ and
$H$ actions are $\Phi_K$ and $\Phi_H= \mathrm{p}\circ \Phi_K$, where
$\mathrm{p}:\kgot^*\to\hgot^*$ is the canonical projection.

\begin{coro}\label{coro:intro-1}
For every $\nu\in\what{H}$, we have :
\begin{equation}\label{eq:coro-restriction}
\Qcal\left(M_{\nu,H}\right)=\sum_{\mu\in \what{K}}N^\mu_{\nu}
\Qcal\left(M_{\mu,K}\right).
\end{equation}
where $M_{\nu,H}=\Phi^{-1}_H(H\cdot\nu)/H$ and
$M_{\mu,K}=\Phi^{-1}_K(K\cdot\mu)/K$ are respectively the symplectic
reductions relative to the $H$ and $K$-actions.
\end{coro}

\medskip

Since $V_\mu^K$ is equal to the $K$-quantization of $K\!\cdot\!\mu$,
the ``quantization commutes with reduction'' Theorem tells us that
$N^\mu_{\nu}= \Qcal( (K\!\cdot\!\mu)_{\nu,H})$ : in particular
$N^\mu_{\nu}\neq 0$ implies that $ \nu\in \mathrm{p}(K\cdot\mu)\
\Longleftrightarrow\ \mu\in K\!\cdot \mathrm{p}^{-1}(\nu)$. Finally
$$
N^\mu_{\nu} \Qcal\left(M_{\mu,K}\right)\neq 0 \ \Longrightarrow\
\mu\in K\!\cdot \mathrm{p}^{-1}(\lambda)\quad \mathrm{and} \quad
\Phi^{-1}_K(\mu)\neq\emptyset.
$$
Theses two conditions imply that
we can restrict the sum of RHS of (\ref{eq:coro-restriction}) to
\begin{equation}\label{eq:intro-mu-phi-nu}
    \mu\in \what{K}\cap \Phi_K\left( K\cdot\Phi^{-1}_H(\nu)\right)
\end{equation}
which is \emph{finite} since $\Phi_H$ is proper.

\medskip

Theorem \ref{theo:intro} and (\ref{eq:functor-produit-formal}) gives
the following extended version of $\textbf{[P3]}$.

\begin{theo}\label{theo:intro-1}
    Let $N$ and $M$ be two pre-quantized
Hamiltonian $K$-manifold where $N$ is \emph{compact} and $M$ is
\emph{proper}. The product $M\times N$ is then proper and we have
have the following equality in $\Rfor(K)$ :
\begin{equation}\label{eq:functor-produit-K}
\qfor_K(M\times N)=\qfor_K(M)\cdot\Qcal_K(N)
\end{equation}
\end{theo}

\medskip

For $\mu,\lambda,\theta\in \what{K}$ we denote
$C^\mu_{\lambda,\theta}=\dim (\hom_H(V_\mu^K,V_\lambda^K\otimes
V_\theta^K))$ the multiplicity of $V_\mu^K$ in the tensor
$V_\lambda^K\otimes V_\theta^K$. Since $V_\lambda^K\otimes
V_\theta^K$ is equal to the quantization of the product
$K\!\cdot\!\lambda\times K\!\cdot\!\theta$, the ``quantization
commutes with reduction'' Theorem tells us that
$C^\mu_{\lambda,\theta}= \Qcal((K\!\cdot\!\lambda\times
K\!\cdot\!\theta)_\mu)$ : in particular $C^\mu_{\lambda,\theta}\neq
0$ implies that $(*)\ \|\lambda\|\leq \|\theta\|+\|\mu\|$.

\begin{coro}\label{coro:intro-2}
In the situation of Theorem \ref{theo:intro-1}, we have for every
$\mu\in\what{K}$ :
\begin{equation}\label{eq:coro-produit}
\Qcal\left((M\times N)_\mu\right)=\sum_{\lambda,\theta\in
\what{K}}C^\mu_{\lambda,\theta}
\Qcal\left(M_\lambda\right)\Qcal\left(N_\theta\right).
\end{equation}
\end{coro}

\medskip

Since $N$ is compact, $\Qcal(N_\theta)\neq 0$ for a $(**)\ \theta\in
\{\mathrm{finite\ set}\}$. Then $(*)$ and $(**)$ show that the sum
in the RHS of (\ref{eq:coro-produit}) is \emph{finite}.

\bigskip

Weistman \cite{Weitsman} studied the formal quantization procedure
in the case of $K=\U(n)$. Using a method of symplectic cutting
\cite{Lerman-cut,Woodward96} he defines for two increasing sequence
of positive integers $r_n,s_n$ a family of cut-spaces $M^{r_n,s_n}$
which are compact. {\em Under the hypothesis that the cut spaces
$M^{r_n,s_n}$ are smooth}, he notes that
$\qfor_K(M)=\lim\limits_{n\to\infty} \Qcal_K(M^{r_n,s_n})$, and he
was then able to show Theorem \ref{theo:intro-1}.

Our proof of Theorem  \ref{theo:intro} uses also a technique of
symplectic cutting but which is valid for any compact Lie group
actions. We have to overpass the difficulties concerning the
non-smoothness of the cut-spaces. For this purpose we introduce
another method of symplectic cutting that uses the wonderfull
compactifications of Concini-Procesi \cite{CP1,CP2}, and we prove an
extension of the ``quantization commutes with reduction'' Theorem to
the \emph{singular} setting.

\bigskip

\textbf{Acknowledgements :}  I am grateful to Michel Brion and
Nicolas Ressayre for enlightening discussions about the wonderfull
compactifications. I thanks also Anton Alekseev for these remarks on
a preliminary version of the paper.

\section{Quantization commutes with reduction}

In this section we precise the definition of the quantization of a
smooth and compact Hamiltonian manifold. We extend the definition to
the case of a \emph{singular} Hamiltonian manifold and we prove a
``quantization commutes with reduction'' Theorem in the singular
setting.

\medskip

In the Kostant-Souriau framework, an Hamiltonian $K$-manifold
$(M,\Omega_M,\Phi_M)$ is pre-quantized if there is an equivariant
Hermitian line bundle $L$ with an invariant Hermitian connection
$\nabla$ such that
\begin{equation}\label{eq:kostant-L}
    \Lcal(X)-\nabla_{X_M}=i\langle\Phi_M,X\rangle\quad \mathrm{and} \quad
    \nabla^2= -i\Omega_M,
\end{equation}
for every $X\in\kgot$. Here $X_M$ is the vector field on $M$ defined
by $X_M(m)=\frac{d}{dt} e^{-tX}m|_{0}$.

$(L,\nabla)$ is also call a Kostant-Souriau line bundle. Remark that
the conditions (\ref{eq:kostant-L}) implies through the equivariant
Bianchi formula the relation
\begin{equation}\label{eq:hamiltonian-action}
    \iota(X_M)\Omega_M= -d\langle\Phi_M,X\rangle,\quad X\in\kgot.
\end{equation}

We will now recall the notions of geometric quantization.

\subsection{Geometric quantization : the compact and smooth case}\label{subsec:quant-smooth}

We suppose here that $(M,\Omega_M,\Phi_M)$ is \textbf{compact} and
is prequantized by a Hermitian line bundle $L$. Choose a
$K$-invariant almost complex structure $J$ on $M$ which is
compatible with $\Omega_M$ in the sense that the symmetric bilinear
form $\Omega_M(\cdot,J\cdot)$ is a Riemannian metric. Let
$\overline{\partial}_L$ be the Dolbeault operator with coefficients
in $L$, and let $\overline{\partial}_L^*$ be its (formal) adjoint.
The \emph{Dolbeault-Dirac operator} on $M$ with coefficients in $L$
is $D_L= \overline{\partial}_L+\overline{\partial}_L^*$, considered
as an operator from $\Acal^{0,\textrm{\tiny even}}(M,L)$ to
$\Acal^{0,\textrm{\tiny odd}}(M,L)$.

\begin{defi}\label{def:quant-compact-lisse}
 The geometric quantization of $(M,\Omega_M,\Phi_M)$ is the element $\Qcal_K(M)\in
R(K)$ which is defined as the equivariant index of the
Dolbeault-Dirac operator $D_L$.
\end{defi}

\begin{rem}\label{rem:dependance-J}
$\bullet$ We can define the Dolbeault-Dirac operator $D_L^{_J}$ for
any invariant almost complex structure $J$. If $J_0$ and $J_1$ are
equivariantly \emph{homotopic} the indices of $D_L^{_{J_0}}$ and
$D_L^{_{J_1}}$ coincides (see \cite{pep-RR}).

$\bullet$ Since the set of \emph{compatible} invariant almost
complex structure on $M$ is path-connected, the element
$\Qcal_K(M)\in R(K)$ does not depend of the choice of $J$.
\end{rem}

\subsection{Geometric quantization : the compact and singular case}\label{subsec:quant-singular}

We are interested to defined the geometric quantization of
\emph{singular} compact Hamiltonian manifolds : here "singular"
means that the manifold is obtain by symplectic reduction.

Let $(N,\Omega_N)$ be a smooth symplectic manifold equipped with an
Hamiltonian action of $K\times H$ : we denote $(\Phi_K,\Phi_H): N\to
\kgot^*\oplus\hgot^*$ the corresponding moment map. We assume that
$N$ is pre-quantized by a $K\times H$-equivariant line bundle $L$
and we suppose that the map $\Phi_H$ is \textbf{proper}. One wants
to define the quantization of the (compact) symplectic quotient
$$
N{/\!\!/}_{0}H:=\Phi_H^{-1}(0)/H.
$$

When $0$ is a regular value of $\Phi_H$, $N{/\!\!/}_{0}H$ is a
compact symplectic \emph{orbifold} equipped with an Hamiltonian
action of $K$ : the corresponding moment map is induced by the
restriction of $\Phi_K$ to $\Phi_H^{-1}(0)$. The symplectic quotient
$N{/\!\!/}_{0}H$ is pre-quantized by the line orbibundle
$$
L_0:= \left(L|_{\Phi_H^{-1}(0)}\right)/H.
$$
Definition \ref{def:quant-compact-lisse} extends to the orbifold
case, so one can still defined the quantization of $N{/\!\!/}_{0}H$
as an element $\Qcal_K(N{/\!\!/}_{0}H)\in R(K)$.

\medskip

Suppose now that $0$ is not a regular value of $\Phi_H$. Let $T_H$
be a maximal torus of $H$, and let $C_H\subset \tgot_H^*$ be a weyl
chamber. Since $\Phi_H$ is proper, the convexity theorem says that
the image of $\Phi_H$ intersects $C_H$ in a closed locally polyhedral 
convex set, that we denoted $\Delta_H(M)$ \cite{L-M-T-W}.

We consider an element $a\in \Delta_H(M)$ which is generic and
sufficiently closed to $0\in \Delta_H(M)$ : we denote $H_a$ the
subgroup of $H$ which stabilizes $a$. When $a\in\Delta_H(M)$  is
generic, one can shows (see \cite{Meinrenken-Sjamaar}) that
$$
N{/\!\!/}_{a}H:=\Phi_H^{-1}(a)/H_a
$$
is a compact $K$-Hamiltonian orbifold, and that
$$
L_a:= \left(L|_{\Phi_H^{-1}(a)}\right)/H_a.
$$
is a $K$-equivariant line orbibundle over $N{/\!\!/}_{a}H$ : we can
then define like in Definition \ref{def:quant-compact-lisse} the
element $\Qcal_K(N{/\!\!/}_{a}H)\in R(K)$ as the equivariant index
of the Dolbeault-Dirac operator on $N{/\!\!/}_{a}H$.

\begin{propdefi}\label{prop:Q-sing}
 The elements $\Qcal_K(N{/\!\!/}_{a}H)\in R(K)$ do not depend of the choice
of the generic element $a\in\Delta_H(M)$, when $a$ is sufficiently
closed to $0$. The common value will be taken as the geometric
quantization of $N{/\!\!/}_{0}H$, and still denoted
$\Qcal_K(N{/\!\!/}_{0}H)$.
\end{propdefi}

\textsc{Proof}. When $N$ is compact and $K=\{e\}$, the proof can be
founded in \cite{Meinrenken-Sjamaar} and in \cite{pep-RR}. The
$\K$-theoretic proof of \cite{pep-RR} extends naturally to our case.
$\Box$

\subsection{Quantization commutes with reduction : the singular case}
\label{subsec:QR-singular}

In section \ref{subsec:quant-singular}, we have defined the
geometric quantization $\Qcal_K(N{/\!\!/}_{0}H)\in R(K)$ of a
compact symplectic reduced space $N{/\!\!/}_{0}H$. We will compute
its $K$-multiplicities like in Theorem \ref{theo:Q-R}.

For every $\mu\in \what{K}$, we consider the coadjoint orbit
$K\cdot\mu\simeq K/K_\mu$ which is pre-quantized by the line bundle
$\C_{[\mu]}\simeq K\times_{K_\mu}\C_\mu$. We consider the
product\footnote{$\overline{K\cdot\mu}$ denotes the coadjoint orbit
with the opposite symplectic form.} $N\times \overline{K\cdot\mu}$
which is an Hamiltonian $K\times H$ manifold which is pre-quantized
by $K\times H$-equivariant line bundle $L\otimes \C_{[\mu]}^{-1}$.
The moment map $N\times \overline{K\cdot\mu}\to
\kgot^*\times\hgot^*, (n,\xi) \mapsto (\Phi_K(n)-\xi,\Phi_H(n))$ is
proper, so the reduced space
$$
\left(N{/\!\!/}_{0}H\right)_\mu:= (N\times\overline{K\cdot\mu})
{/\!\!/}_{(0,0)}K\times H
$$
is compact. Following Proposition \ref{prop:Q-sing}, we can then
define its quantization \break
$\Qcal(\left(N{/\!\!/}_{0}H\right)_\mu)\in\Z$. The main result of
this section is the
\begin{theo}\label{theo:Q-R-singular}
We have the following equality in $R(K)$ :
\begin{equation}\label{eq:Q-R-singular}
    \Qcal_K(N{/\!\!/}_{0}H)=\sum_{\mu\in
\what{K}}\Qcal(\left(N{/\!\!/}_{0}H\right)_\mu) V_\mu^K.
\end{equation}
\end{theo}

\textsc{Proof}. The proof will occupied the remaining of this
section. The starting point is to state another definition of the
geometric quantization of a symplectic reduced space which uses the
Atiyah-Singer's theory of transversally elliptic operators.

\subsubsection{Transversally elliptic symbols}\label{subsec:transversally}

Here we give the basic definitions from the theory of transversally
elliptic symbols (or operators) defined by Atiyah-Singer in
\cite{Atiyah74}. For an axiomatic treatment of the index morphism
see Berline-Vergne \cite{B-V.inventiones.96.1,B-V.inventiones.96.2}
and for a short introduction see \cite{pep-RR}.

Let $\Xcal$ be a {\em compact} $K_1\times K_2$-manifold. Let $p:\T
\Xcal\to \Xcal$ be the projection, and let $(-,-)_\Xcal$ be a
$K_1\times K_2$-invariant Riemannian metric. If $E^{0},E^{1}$ are
$K_1\times K_2$-equivariant vector bundles over $\Xcal$, a
$K_1\times K_2$-equivariant morphism $\sigma \in \Gamma(\T
\Xcal,\hom(p^{*}E^{0},p^{*}E^{1}))$ is called a {\em symbol}. The
subset of all $(x,v)\in \T \Xcal$ where $\sigma(x,v): E^{0}_{x}\to
E^{1}_{x}$ is not invertible is called the {\em characteristic set}
of $\sigma$, and is denoted by $\Char(\sigma)$.

Let $\T_{K_2}\Xcal$ be the following subset of $\T \Xcal$ :
$$
   \T_{K_2}\Xcal\ = \left\{(x,v)\in \T M,\ (v,X_{M}(x))_{_{M}}=0 \quad {\rm for\ all}\
   X\in\kgot_2 \right\} .
$$

A symbol $\sigma$ is {\em elliptic} if $\sigma$ is invertible
outside a compact subset of $\T \Xcal$ (i.e. $\Char(\sigma)$ is
compact), and is $K_2$-{\em transversally elliptic} if the
restriction of $\sigma$ to $\T_{K_2}\Xcal$ is invertible outside a
compact subset of $\T_{K_2}\Xcal$ (i.e. $\Char(\sigma)\cap
\T_{K_2}\Xcal$ is compact). An elliptic symbol $\sigma$ defines an
element in the equivariant $\K$-theory of $\T\Xcal$ with compact
support, which is denoted by $\K_{K_1\times K_2}(\T \Xcal)$, and the
index of $\sigma$ is a virtual finite dimensional representation of
$K_1\times K_2$
\cite{Atiyah-Segal68,Atiyah-Singer-1,Atiyah-Singer-2,Atiyah-Singer-3}.

A $K_2$-{\em transversally elliptic} symbol $\sigma$ defines an
element of $\K_{K_1\times K_2}(\T_{K_2}M)$, and the index of
$\sigma$ is defined as a trace class virtual representation of
$K_1\times K_2$ (see \cite{Atiyah74} for the analytic index and
\cite{B-V.inventiones.96.1,B-V.inventiones.96.2} for the
cohomological one) : in fact $\indice^\Xcal(\sigma)$ belongs to the
tensor product $R(K_1)\what{\otimes} \Rfor(K_2)$.

Remark that any elliptic symbol of $\T \Xcal$ is $K_2$-transversally
elliptic, hence we have a restriction map $\K_{K_1\times K_2}(\T
\Xcal)\to \K_{K_1\times K_2}(\T_{K_2}\Xcal)$, and a commutative
diagram
\begin{equation}\label{indice.generalise}
\xymatrix{ \K_{K_1\times K_2}(\T \Xcal) \ar[r]\ar[d]_{\indice^\Xcal}
&
\K_{K_1\times K_2}(\T_{K_2}\Xcal)\ar[d]^{\indice^\Xcal}\\
R(K_1)\otimes R(K_2)\ar[r] &  R(K_1)\what{\otimes}\ \Rfor(K_2)\ .
   }
\end{equation}

\medskip

Using the {\em excision property}, one can easily show that the
index map $\indice^\Ucal: \K_{K_1\times K_2}(\T_{K_2}\Ucal)\to
R(K_1)\what{\otimes} \Rfor(K_2)$ is still defined when $\Ucal$ is a
$K_1\times K_2$-invariant relatively compact open subset of a
$K_1\times K_2$-manifold (see \cite{pep-RR}[section 3.1]).

\subsubsection{Quantization of singular space : second definition}

Let $(\Xcal,\Omega_\Xcal)$ be an Hamiltonian $K_1\times
K_2$-manifold pre-quantized by a $K_1\times K_2$-equivariant line
bundle $L$. The moment map $\Phi_2: \Xcal\to\kgot_2^*$ relative to
the $K_2$-action is supposed to be \textbf{proper.} Take a
compatible $K_1\times K_2$-invariant almost complex structure on
$\Xcal$. We choose a $K_1\times K_2$-invariant Hermitian metric
$\|v\|^2$ on the tangent bundle $\T \Xcal$, and we identify the
cotangent bundle $\T^*\Xcal$ with $\T \Xcal$. For $(x,v)\in\T\Xcal$,
the principal symbol of the Dolbeault-Dirac operator
$\overline{\partial}_L+\overline{\partial}_L^*$ is the clifford
multiplication $\clif_\Xcal(v)$ on the complex vector bundle
$\Lambda^\bullet\T_x \Xcal\otimes L_x$. It is invertible for $v\neq
0$, since $\clif_\Xcal(v)^2=-\|v\|^2$.

When $\Xcal$ is compact, the symbol $\clif_\Xcal$ is elliptic and
then defines an element of the equivariant \textbf{K}-group of
$\T\Xcal$. The topological index of $\clif_\Xcal\in\K_{K_1\times
K_2}(\T\Xcal)$ is equal to the analytical index of the
Dolbeault-Dirac operator
$\overline{\partial}_L+\overline{\partial}_L^*$ :
\begin{equation}\label{eq:index-analytique-topologique}
    \Qcal_{K_1\times K_2}(\Xcal)=\indice^\Xcal(\clif_\Xcal)\quad
    \mathrm{in}\quad R(K_1)\otimes R(K_2).
\end{equation}

\medskip

When $\Xcal$ is not compact the topological index of $\clif_\Xcal$
is not defined. In order to give a topological definition of
$\Qcal_{K_1}(\Xcal{/\!\!/}_{0}K_2)$, we will deform the symbol
$\clif_\Xcal$ as follows. Consider the identification
$\kgot_2^*\simeq\kgot_2$ defined by a $K_2$-invariant scalar product
on the Lie algebra $\kgot_2$. From now on the moment map $\Phi_2$
will take values in $\kgot_2$, and we define the vectors field on
$\Xcal$
\begin{equation}\label{eq-kappa}
    \kappa_x= (\Phi_2(x))_M(x), \quad x\in\Xcal.
\end{equation}
We consider now the symbol
$$
\clif_\Xcal^\kappa(v)=\clif(v-\kappa_x), \quad v\in \T_x\Xcal.
$$
Note that $\clif_\Xcal^\kappa(v)$ is invertible except if
$v=\kappa_x$. If furthermore $v$ belongs to the subset $\T_{K_2} M$
of tangent vectors orthogonal to the $K_2$-orbits, then $v=0$ and
$\kappa_x=0$.  Indeed $\kappa_x$ is tangent to $K_2\cdot x$ while
$v$ is orthogonal.

Since $\kappa$ is the Hamiltonian vector field of the function
$\|\Phi_2\|^2$, the set of zeros of $\kappa$ coincides with the set
$\Cr(\|\Phi_2\|^2)$ of critical points of $\|\Phi_2\|^2$.

Let $\Ucal\subset\Xcal$ be a $K_1\times K_2$-invariant open subset
which is \emph{relatively compact}. If the border $\partial\Ucal$
does not intersect $\Cr(\|\Phi_2\|^2)$, then the restriction
$c^\kappa_\Xcal|_\Ucal$ defines a class in $\K_{K_1\times
K_2}(\T_{K_2}\Ucal)$ since
$$
\Char(\clif_\Xcal^\kappa|_\Ucal)\cap \T_{K_2}\Ucal \simeq
\Cr(\|\Phi_2\|^2)\cap\Ucal
$$
is compact. In this situation the index of $\clif_\kappa|_\Ucal$ is
defined as an element $\indice^\Ucal(\clif_\kappa|_\Ucal)\in
R(K_1)\what{\otimes} \Rfor(K_2)$.

\begin{theo}\label{theo:quant-sing-bis}
    The $K_2$-invariant part of $\indice^\Ucal(\clif_\Xcal^\kappa|_\Ucal)$ is equal
    to :
    \begin{itemize}
      \item $\Qcal_{K_1}(\Xcal{/\!\!/}_{0}K_2)$ when
      $\Phi_2^{-1}(0)\subset \Ucal$,
      \item $0$ in the other case.
    \end{itemize}
\end{theo}

\textsc{Proof}. When $K_1=\{e\}$, the proof is done in \cite{pep-RR}
(see section 7). This proof works equally well in the general case.

\begin{rem}\label{rem:X-compact}
    If $\Xcal$ is compact we can take $\Ucal=\Xcal$ in the last
    Theorem. In this case the symbols $\clif_\Xcal^\kappa$ and $\clif_\Xcal$
    defines the same class in $\K_{K_1\times K_2}(\T\Xcal)$ so they
    have the same index. Theorem \ref{theo:quant-sing-bis}
    corresponds then to the traditional ``quantization commutes with
    reduction'' phenomenon : $[\Qcal_{K_1\times
K_2}(\Xcal)]^{K_2}=\Qcal_{K_1}(\Xcal{/\!\!/}_{0}K_2)$.
\end{rem}

>From now one  we will work with this topological definiton for the
geometric quantization of the reduced $K_1$-Hamiltonian manifold
$\Xcal{/\!\!/}_{0}K_2$ (which is possibly singular):
$\Qcal_{K_1}(\Xcal{/\!\!/}_{0}K_2)=
[\indice^\Ucal(\clif_\Xcal^\kappa|_\Ucal)]^{K_2}$ where $\Ucal$ is
any relatively compact neighborhood of $\Phi_2^{-1}(0)$ such that
$\partial\Ucal\cap \Cr(\|\Phi_2\|^2)=\emptyset$. The functorial
properties still holds in this singular setting. In particular :

\medskip

$\textbf{[P2]}$  If $H\subset K_1$ is a connected lie subgroup, then
the restriction of $\Qcal_{K_1}(\Xcal{/\!\!/}_{0}K_2)$ to $H$ is
equal to $\Qcal_{H}(\Xcal{/\!\!/}_{0}K_2)$.


\subsubsection{Proof of theorem \ref{theo:Q-R-singular}}

We come back in the situation of sections
\ref{subsec:quant-singular} and \ref{subsec:QR-singular}~.

First we apply Theorem \ref{theo:quant-sing-bis} to $\Xcal=N$,
$K_1=K$ and $K_2=H$. (\ref{eq:Q-R-singular}) is trivially true when
$0\notin \image(\Phi_H)$. So we suppose now that $0\in
\image(\Phi_H)$, and we consider an $K\times H$-invariant open
subset $\Ucal\subset N$ which is \emph{relatively compact} and such
that
$$
\Phi_H^{-1}(0)\subset \Ucal\quad \mathrm{and} \quad
\partial\Ucal\cap \Cr(\|\Phi_H\|^2)=\emptyset.
$$
We have $\Qcal_{K}(N{/\!\!/}_{0}H)=[\indice^\Ucal
(\clif_N^{\kappa^{_H}}|_\Ucal)]^H$ and one want to compute its $K$-
multiplicities $m_\mu, \mu\in \what{K}$. Here $\kappa^{_H}$ is the
vectors field on $N$ associated to the moment map $\Phi_H$ (see
(\ref{eq-kappa})).

Take $\mu\in \what{K}$. We denote $\clif_{-\mu}$ the principal
symbol of the Dolbeault-Dirac operator on $\overline{K\cdot\mu}$
with values in the line bundle $\C_{[-\mu]}$ : we have
$\indice^{K\cdot\mu}(\clif_{-\mu})=(V_\mu^{K})^*$.

We know then that the multiplicity  of
$[\indice^\Ucal(\clif_N^{\kappa^{_H}}|_\Ucal)]^H$ relatively to
$V_\mu^{K}$ is equal to
\begin{equation}\label{eq:m-mu-1}
    m_\mu:=\Big[\indice^{\Vcal}
\left(\clif_N^{\kappa^{_H}}|_\Ucal\odot\clif_{-\mu}\right)\Big]^{K\times
H}
\end{equation}
with $\Vcal=\Ucal\times K\cdot\mu$. This identity is due to the fact
that we have a "multiplication"
\begin{eqnarray*}\label{eq:produit.transversal}
\K_{K\times H}(\T_{H}\Ucal)\times \K_{K}(\T (K\!\cdot\!\mu))
&\longrightarrow &
\K_{K\times H}(\T_{K\times H}(\Ucal\times K\!\cdot\!\mu) )\\
(\sigma_{1},\sigma_{2})&\longmapsto &  \sigma_{1}\odot \sigma_{2} \
.\nonumber
\end{eqnarray*}
so that $\indice^{\Ucal\times K\cdot\mu}(\sigma_{1}\odot\sigma_{2})=
\indice^{\Ucal}(\sigma_{1})\cdot\indice^{K\cdot\mu}(\sigma_{2})$ in
$R^{-\infty}(K\times H)$. See \cite{Atiyah74}.

\bigskip

Consider now the case where  $\Xcal=N\times \overline{K\cdot\mu} $,
$K_1=\{e\}$ and $K_2=K\times H$. After Theorem
\ref{theo:quant-sing-bis}, we know that
\begin{equation}\label{eq:m-mu-2}
   \Qcal(\left(N{/\!\!/}_{0}H\right)_\mu)
   = \Big[\indice^\Vcal(\clif_\Xcal^\kappa|_\Vcal)\Big]^{K\times H},
\end{equation}
where $\kappa$ is the vector field on $N\times \overline{K\cdot\mu}$
associated to the moment map
\begin{eqnarray*}
  \Phi : N\times \overline{K\cdot\mu} &\longrightarrow & \kgot^*\times\hgot^* \\
   (x,\xi) &\longmapsto& (\Phi_K(x)-\xi, \Phi_H(n))
\end{eqnarray*}
Note that $\Vcal=\Ucal\times K\cdot\mu$ is a neighborhood of
$\Phi^{-1}(0)\subset (\Phi_H)^{-1}(0)$.

\medskip

Our aim now is to prove that the quantities (\ref{eq:m-mu-1}) and
(\ref{eq:m-mu-2}) are equal.

\medskip

Since the definition of $\kappa$ needs the choice of an invariant
scalar product on the Lie algebra $\kgot\times\hgot$, we will
precise its definition. Let $\|\cdot\|_K$ and $\|\cdot\|_H$ be two
invariant Euclidean norm respectively on $\kgot$ and $\hgot$. For
any $r>0$ we consider on $\kgot\times\hgot$ the invariant Euclidean
norm $\|(X,Y)\|^2_r=r^2\|X\|^2_K+  \|Y\|_H^2$.

Let $\kappa^{_K}$ be the vector field on $N\times
\overline{K\cdot\mu}$ associated to the map $N\times
\overline{K\cdot\mu} \to\kgot^*, (x,\xi) \mapsto \Phi_K(x)-\xi$, and
where the identification $\kgot\simeq \kgot^*$ is made through the
Euclidean norm $\|\cdot\|_K$ (see (\ref{eq-kappa})). For $(x,\xi)\in
N\times \overline{K\cdot\mu}$, we have the decomposition
$$
\kappa^{_K}(x,\xi)=(\kappa_1(x,\xi),\kappa_2(x,\xi))\ \in\  \T_x
N\times \T_\xi K\!\cdot\!\mu.
$$

Let $\kappa^{_H}$ be the vector field on $N\times
\overline{K\cdot\mu}$ associated to the map $N\times
\overline{K\cdot\mu} \to\kgot^*, (x,\xi) \mapsto \Phi_H(x)$, and
where the identification $\hgot\simeq \hgot^*$ is made through the
Euclidean norm $\|\cdot\|_H$. For $(x,\xi)\in N\times
\overline{K\cdot\mu}$, we have the decomposition
$$
\kappa^{_H}(x,\xi)=(\kappa^{_H}(x),0)\ \in\  \T_x N\times \T_\xi
K\!\cdot\!\mu.
$$

For any $r> 0$, we denote by $\kappa_r$ the vector field on $N\times
\overline{K\cdot\mu}$ associated to the map $\Phi$, and where the
identification $\kgot\times\hgot\simeq \kgot^*\times\hgot^*$ is made
through the Euclidean norm $\|\cdot\|_r$. We have then
\begin{eqnarray*}
  \kappa_r &=& \kappa^{_H} +r\, \kappa^{_K} \\
         &=& (\kappa^{_H} +r\,\kappa_1,r\,\kappa_2)
\end{eqnarray*}

Now we can precise (\ref{eq:m-mu-2}). Take an invariant relatively
compact neighborhood $\Ucal$ of $\Phi_H^{-1}(0)$ such that $\partial
\Ucal \cap\{\mathrm{zeros\ of}\ \kappa^{_H} \}=\emptyset$. With the
help of a invariant Riemannian metric on $\Xcal$ we define
$$
\varepsilon_H=\inf_{x\in \partial\Ucal}\|\kappa^{_H}(x)\|\ > 0\quad
\mathrm{and} \quad \varepsilon_K=\sup_{(x,\xi)\in
\partial\Ucal\times K\!\cdot\!\mu }\|\kappa_1(x,\xi)\|.
$$
Note that for any $0\leq r<\frac{\varepsilon_H}{\varepsilon_K}$, we
have $\partial \Ucal \times K\!\cdot\!\mu\cap\{\mathrm{zeros\ of}\
\kappa^{_H}+r\kappa_1 \}=\emptyset$, and then $\partial \Vcal
\cap\{\mathrm{zeros\ of}\ \kappa_r \}=\emptyset$ for the
neighborhood $\Vcal:=\Ucal\times K\cdot\mu$ of $\Phi^{-1}(0)$. We
can then use Theorem \ref{theo:quant-sing-bis} : for
$0<r<\frac{\varepsilon_H}{\varepsilon_K}$ we have
$$
\Qcal(\left(N{/\!\!/}_{0}H\right)_\mu)=
\Big[\indice^\Vcal(\clif_\Xcal^{\kappa_r}|_\Vcal)\Big]^{K\times H}
$$
We are now close to the end of the proof. Let us compare the symbols
$\clif_\Xcal^{\kappa_r}|_\Vcal$ and
$\clif_N^{\kappa^{_H}}|_\Ucal\odot\clif_{-\mu}$ in $\K_{K\times
H}(\T_{K\times H}(\Ucal\times K\!\cdot\!\mu))$. First one sees that
the symbols $\clif_\Xcal$ is equal to the product
$\clif_N\odot\clif_{-\mu}$ hence the symbols
$\clif_N^{\kappa^{_H}}|_\Ucal\odot\clif_{-\mu}$ is equal to
$\clif_\Xcal^{\kappa_r}|_\Vcal$ when $r=0$. Since for
$r<\frac{\varepsilon_H}{\varepsilon_K}$ the path $s\in[0,r]\to
\clif_\Xcal^{\kappa_s}|_\Vcal$ defines an homotopy of $K\times
H$-transversally elliptic symbols on $\Vcal$, we get
$$
\indice^\Vcal(\clif_\Xcal^{\kappa_r}|_\Vcal)=
\indice^\Vcal(\clif_N^{\kappa^{_H}}|_\Ucal\odot\clif_{-\mu})
$$
and then $m_\mu=\Qcal(\left(N{/\!\!/}_{0}H\right)_\mu)$. $\Box$

\section{Wonderful compactifications and symplectic cutting}

In this section we use the wonderful compactifications of
Concini-Procesi \cite{CP1,CP2} to perform symplectic cutting.

\subsection{Wonderful compactifications : definitions}

We study here the wonderfull compactifications from the Hamiltonian
point of view.

We consider a compact connected Lie group $K$ and its
complexification $K_\C$. Let $T$ be a maximal torus of $K$, and let
$W:=N(T)/T$ be the Weyl group. Consider $\tgot^*$, the dual of the
Lie algebra of $T$, with the lattice $\wedge^*$ of real weights. Let
$C_K\subset \tgot^*$ be a Weyl chamber and let $\what{K}:=
\wedge^*\cap C_K$ be the set dominants weights.


\begin{defi}\label{defi:polytope-adap}
A polytope $P$ in $\tgot^*$ is $K$-adapted if :
\begin{itemize}
  \item[i)] the vertices of $P$ are regular elements of $\wedge^*$,
  \item[ii)] $P$ is $W$-invariant,
  \item[iii)] $P$ is Delzant.
\end{itemize}
\end{defi}

\textsc{Example} : When $K$ as a trivial center, the convex hull of
$W\cdot\mu$ is a $K$-adapted polytope for any regular dominant
weight $\mu$.

\begin{prop}
There exists $K$-adapted polytopes in $\tgot^*$.
\end{prop}

\textsc{Proof.} Let us use the dictionary between polytopes and
projective fan \cite{Oda88}. Conditions $ii)$ and $iii)$ of Definition
\ref{defi:polytope-adap} means that we are looking after a smooth
projective $W$-invariant fan $\Fcal$ in $\tgot$. Condition $i)$
means that each cone of $\Fcal$ of maximal dimension should not be
fixed by any element of $W\setminus\{{\rm Id}\}$. For a proof of the
existence of search fan, see \cite{Brylinski79,CT-H-S}. In
particular condition $(*)$ in Proposition 2 of \cite{CT-H-S} implies
$i)$. $\Box$

\bigskip

In the rest of this section, we consider a $K$-adapted polytope $P$.
Let $[P]_1$ be the union of all the closed facet of dimension $1$ :
we label the elements of $[P]_1\cap\what{K}$ by
$\{\lambda_1,\cdots,\lambda_N\}$. Since $P$ is $K$-adapted, when
$\lambda_{i}$ is a vertex of $P$ there exits
$\alpha_{j_1},\cdots,\alpha_{j_r}$ belonging to $[P]_1\cap\wedge^*=
W\cdot\{\lambda_1,\cdots,\lambda_N\}$ such that
$\alpha_{j_1}-\lambda_{i},\cdots,\alpha_{j_r}-\lambda_{i}$ is a
basis of the lattice $\wedge^*$.

Let $V_{\lambda_i}$ be an irreducible representation of $K$ with
eighest weight $\lambda_i$: these representations extend canonically
to the complexification $K_\C$. We denote $\rho: K_\C\to \Pi_{i=1}^N
\GL(V_{\lambda_i})$ the representation of $K_\C$ on $V:=\oplus_{i=1}^N
V_{\lambda_i}$. We consider the vector space
$$
E=\oplus_{i=1}^N \End(V_{\lambda_i})
$$
equipped with the action of $K_\C\times K_\C$ given by :
$(g_1,g_2)\cdot f=\rho(g_1)\circ f\circ \rho(g_2)^{-1}$. Let
$\Pbb(E)$ the projective space associated to $E$ : it is equipped
with an algebraic action of the reductive group $K_\C\times K_\C$.
We consider the map $g\mapsto [\rho(g)]$ from $K_\C$ into $\Pbb(E)$,
that we denote $\bar{\rho}$.

\begin{lem}\label{lem:properties-rho}
The map $\bar{\rho}: K_\C\to \Pbb(E)$ is an embedding.

\end{lem}

\textsc{Proof}. Let $g\in K_\C$ such that $\bar{\rho}(g)= [Id]$ :
there exist $a\in \C^*$ such that $\rho(g)=a Id$. The Cartan
decomposition gives
\begin{equation}\label{eq:cartan-rho}
\rho(k)=\frac{a}{|a|}Id\quad \mathrm{and}\quad\rho(e^{iX})= |a| Id
\end{equation}
for $g=ke^{iX}$ with $k\in K$ and $X\in\kgot$. Since there exist
$Y,Y'\in \tgot$ and $u,u'\in K$ such that $k=u e^Y u^{-1}$ and
$X=u'\cdot Y'$, (\ref{eq:cartan-rho}) gives
\begin{equation}\label{eq:cartan-rho-bis}
    \rho(e^Y)=\frac{a}{|a|}Id\quad \textrm{and} \quad
\rho(e^{iY'})=|a|Id.
\end{equation}

Since each element of $[P]_1\cap\wedge^*=
W\cdot\{\lambda_1,\cdots,\lambda_N\}$ is a weight for the action of
$T_\C$ on $\oplus_{i=1}^N V_{\lambda_i}$, (\ref{eq:cartan-rho-bis})
implies that for every $\alpha\in [P]_1\cap\wedge^*$ we have
\begin{equation}\label{eq:exp-y}
    e^{i\langle \alpha,Y\rangle}=\frac{a}{|a|}\quad \mathrm{and}\quad
e^{-\langle\alpha,Y'\rangle}=|a|.
\end{equation}
Since there exists $\alpha_{j_0},\cdots,\alpha_{j_r}\in
[P]_1\cap\wedge^*$ such that
$\alpha_{j_1}-\alpha_{j_0},\cdots,\alpha_{j_r}-\alpha_{j_0}$ is a
basis of the lattice $\wedge^*$, (\ref{eq:exp-y}) implies that
$Y'=0$ and that $Y\in \ker(Z\in\tgot\to e^Z)$. We have proved that
$g=e$. $\Box$


\medskip

Let $T_\C\subset K_\C$ the complexification of the (compact) torus
$T\subset K$.

\begin{defi}\label{def:cp}
    Let $\Xcal_P$ be the Zarisky closure of $\bar{\rho}(K_\C)$ in $\Pbb(E)$ and let
$\Ycal_P\subset \Xcal_P$ be the Zarisky closure of
$\bar{\rho}(T_\C)$ in $\Pbb(E)$.
\end{defi}

Since $\bar{\rho}(K_\C)=K_\C\times K_\C\cdot [Id]$ and
$\bar{\rho}(T_\C)=T_\C\times T_\C\cdot [Id]$ are orbits of algebraic
group actions their Zariski closures coincide with their closures
for the Euclidean topology.

\begin{theo}\label{th:cp-lisse}
    The varities $\Xcal_P$ and $\Ycal_P$ are smooth.
\end{theo}

The proof will be done in the next section

\subsection{Smoothness of $\Xcal_P$ and $\Ycal_P$}

Let $E$ be a complex vector space equipped with a linear action of a
reductive group $G$. Let $\Zcal\subset \Pbb(E)$ be a projective
variety which is $G$-stable. We have the classical fact
\begin{lem}\label{lem:variety-smooth}
    \begin{itemize}
      \item $\Zcal$ possess closed $G$-orbits.
      \item $\Zcal$ is smooth if $\Zcal$ is smooth near its closed $G$-orbits.
    \end{itemize}
\end{lem}

\textsc{Proof}. Let $z_0\in \Zcal$ and consider the Zariski closure
$\overline{G\cdot z_0}\subset \Zcal$. If $G\cdot z_0$ is not closed,
we take $z_1\in \overline{G\cdot z_0}\setminus G\cdot z_0$ : we have
$\dim G\cdot z_1< \dim G\cdot z_0$. By induction we find a sequence
$z_1,\cdots, z_p$ with $z_{k+1}\in \overline{G\cdot z_k}\setminus
G\cdot z_k$ for $k<p$ and $G\cdot z_p$ closed. For the second point,
we have just to note that if $\Zcal$ is singular, the subvariety
$\Zcal^{sing}\subset\Zcal$ of singular points is $G$-stable and then
contains a closed $G$-orbits. $\Box$

\medskip

We are interested here respectively in

\begin{itemize}
  \item the $K_\C\times K_\C$-variety $\Xcal_P\subset\Pbb(E)\subset\Pbb(\End(V))$
  \item the $T_\C\times T_\C$-variety $\Ycal_P\subset\Pbb(E)$.
\end{itemize}
Since the diagonal $Z_\C=\{(t,t)| t\in T_\C\}$
 stabilizes $[Id]$, its action on $\Ycal_P$ is trivial. Hence we will restrict
ourself to the action of $T_\C\times T_\C/Z_\C\simeq T_\C$ on
$\Ycal_P$ : for $t\in T_\C$ and $[y]\in\Ycal_P$ we take
$t\cdot[y]=[\rho(t)\circ y]$.

\subsubsection{The case of $\Ycal_P$}

We apply Lemma \ref{lem:variety-smooth} to the $T_\C$-variety
$\Ycal_P=\overline{T_\C\cdot[Id]}$ in $\Pbb(E)$. Let $\{\alpha_j,
j\in J\}$ be the $T_\C$ weights on $(\oplus_{i=1}^N
V_{\lambda_i},\rho)$ counted with their multiplicity. Their exists a
orthonormal basis $\{v_j,j\in J\}$ of $\oplus_{i=1}^N V_{\lambda_i}$
such that $Id=\sum_{j\in J}v_j\otimes v_j^*$ and
\begin{equation}\label{eq:rho-T}
    \rho(e^Z)=\sum_{j\in J} e^{i\langle\alpha_j,Z\rangle}v_j\otimes
    v_j^*,\quad Z\ \in \tgot_\C.
\end{equation}
So the action of $e^Z\in T_\C$ on $[Id]\in\Pbb(E)$ is
$e^Z\cdot[Id]=\left[\sum_{j\in J}
e^{i\langle\alpha_j,Z\rangle}v_j\otimes
    v_j^*\right]$.
We introduce a subset $J'$ of $J$ such that for every $j\in J$ there
exists a \emph{unique} $j'\in J'$ such that $\alpha_j=\alpha_{j'}$.
So the variety $\Ycal_P$ belongs to $\Pbb(E')$ where
$E'=\oplus_{j'\in J'} \C m_{j'}$ with
$m_{j'}=\sum_{j,\alpha_j=\alpha_{j'}} v_j\otimes v_j^*$. The closed
$T_\C$ orbits in $\Pbb(E')$ are $[m_{j'}]$, $j'\in J$.

\begin{lem}\label{lem:Y-closed-orbit}
$[m_{j_o}]\in \Ycal_P$ if and only if $\alpha_{j_o}$ is a vertex of
the polytope $P$.
\end{lem}

\textsc{Proof}. If $\alpha_{j_o}$ is a vertex of $P$, there exists
$X\in \tgot$ such that
$\langle\alpha_{j_o},X\rangle>\langle\alpha_{j},X\rangle$ whenever
$\alpha_{j_o}\neq\alpha_{j}$. Hence $e^{-isX}\cdot [Id]$ tends to
$[m_{j_o}]$ when $s\to+\infty$. If $\alpha_{j_o}$ is not a vertex of
$P$, there exist $L\subset J'\setminus \{j_o\}$ such that
$\alpha_{j_o}=\sum_{l\in L}a_l\alpha_l$ with $0<a_l<1$ and $\sum_l
a_l=1$. So $\Ycal_P$ belongs to the closed subset defined by
$$
[\sum_{j'\in J'}\delta_{j'} m_{j'}]\in \Pbb(E')\quad \mathrm{such\
that}\quad \Pi_{l\in L}|\delta_l|^{a_l}=|\delta_{j_o}|.
$$
Hence $[m_{j_o}]\notin \Ycal_P$. $\Box$

\begin{rem}\label{rem:vertex}
    When $\alpha_{j}$ is a vertex of
the polytope $P$, the multiplicity of $\alpha_{j}$ in
$\oplus_{i=1}^N V_{\lambda_i}$ is equal to one, so $m_{j}=
v_j\otimes v_j^*$.
\end{rem}

\medskip

Consider now a vertex $\alpha_{j_o}$ of $P$ (for $j_o\in J'$). We
consider the open subset $\Vcal\subset\Pbb(E')$ defined by
$[\sum_{j'\in J'}\delta_{j'} m_{j'}]\in \Vcal \Leftrightarrow
\delta_{j_o}\neq 0$, and the diffeomorphism $\psi:\Vcal\to
\C^{J'\setminus \{j_o\}}$, $[\sum_{j'\in J'}\delta_{j'}
m_{j'}]\mapsto (\frac{\delta_{j'}}{\delta_{j_o}})_{j'\neq j_o}$. The
map $\psi$ realizes an  algebraic diffeomorphism between
$\Ycal_P\cap\Vcal$ and the affine subvariety
$$
\Zcal:=\overline{\{(t^{\alpha_{j'}-\alpha_{j_o}})_{j'\neq j_o}\ | \
t\in T_\C\}}\subset\C^{J'\setminus \{j_o\}}.
$$

The set of weights $\alpha_j,j\in J$ contains all the lattice points
that belongs to the one dimensional faces of $P$. Since the polytope
$P$ is $K$-adapted, there exists a subset $L_{j_o}\subset J'$ such
that $\alpha_l-\alpha_{j_o}, l\in L_{j_o}$ is a $\Z$-basis of the
group of weights $\wedge^*$. And for every $j'\neq j_o$ we have
\begin{equation}\label{eq:alpha-base}
\alpha_{j'}-\alpha_{j_o}=\sum_{l\in L_{j_o}} n^l_{j'}
(\alpha_l-\alpha_{j_o})\quad \textrm{with} \quad n^l_{j'}\in \N.
\end{equation}
We define on $\C^{L_{j_o}}$ the monomials $P_{j'}(Z)=\Pi_{l\in
L_{j_o}}(Z_l)^{n^l_{j'}}$. Note that $P_{j'}(Z)=Z_{l}$ when $j'=l\in
L_{j_o}$. Now it is not difficult to see that the map
\begin{eqnarray*}
  \C^{L_{j_o}} &\longrightarrow & \C^{J'\setminus \{j_o\}} \\
  Z &\longmapsto & \left(P_{j'}(Z)\right)_{j'\neq j_o}
\end{eqnarray*}
realizes an algebraic diffeomorphism between $\C^{L_{j_o}}$ and
$\Zcal$.

\bigskip

Finally we have shown that $\Ycal_P$ is smooth near $[m_{j_o}]$ :
hence $\Ycal_P$ is a smooth subvariety of $\Pbb(E)$. Since $T_\C$
acts on $\Ycal_P$ with a dense orbit, $\Ycal_P$ is a smooth
projective toric variety.


\subsubsection{The case of $\Xcal_P$}

Let $E:=\otimes_{i=1}^N \End(V_{\lambda_i})$. The closed $K_\C\times
K_\C$-orbit in $\Pbb(E)$ are those passing through
$[v_{\lambda_i}\otimes v_{\lambda_i}^*]$ where $v_{\lambda_i}\in
V_{\lambda_i}$ is a highest weight vector (that we take of
norm $1$ for a $K$-invariant hermitian structure).

\begin{lem}\label{lem:X-closed-orbit}
$[v_{\lambda_i}\otimes v_{\lambda_i}^*]\in \Xcal_P$ if and only if
$\lambda_i$ is a vertex of the polytope $P$.
\end{lem}

\textsc{Proof}. If $\lambda_i$ is a vertex of $P$, we have proved in
Lemma \ref{lem:Y-closed-orbit} that $[v_{\lambda_i}\otimes
v_{\lambda_i}^*]$ belongs to  $\Ycal_P$ and so belongs to $\Xcal_P$.
We prove the converse in Corollary \ref{coro:vertex-P}. $\Box$

\medskip

For the remaining of this section we consider a vertex
$\lambda_{i_o}\in \what{K}$ of the polytope $P$. Let $B^+, B^-$ the
Borel subgroups fixing respectively the elements
$[v_{\lambda_{i_o}}]\in\Pbb(V_{\lambda_i})$ and
$[v_{\lambda_{i_o}}^*]\in\Pbb(V_{\lambda_i}^*)$. Consider also the
unipotent subgroups $N^\pm\subset B^\pm$ fixing respectively the
elements $v_{\lambda_i}\in V_{\lambda_i}$ and $v_{\lambda_i}^*\in
V_{\lambda_i}^*$.

We consider the open subset $\Vcal_{\textrm{End}} \subset \Pbb(E)$
of elements $[f]$ such that \break $\langle
v_{\lambda_{i_o}}^*,f(v_{\lambda_{i_o}})\rangle\neq 0$ :
$\Vcal_{\textrm{End}}$ is $B^-\times B^+$ stable. Consider the open
subset $\Vcal \subset \Pbb(V_{\lambda_{i_o}})$ and $\Vcal^* \subset
\Pbb(V_{\lambda_{i_o}}^*)$) defined by :
\begin{itemize}
  \item $[v]\in\Vcal\ \Leftrightarrow\ \langle v_{\lambda_{i_o}}^*, v\rangle\neq
  0$ : $\Vcal$ is $B^-$ stable,
  \item $[\xi]\in\Vcal^*\ \Leftrightarrow\  \langle \xi,v_{\lambda_{i_o}}\rangle\neq
  0$ : $\Vcal^*$ is $B^+$ stable.
\end{itemize}
We consider now the rational maps $l:\Pbb(E)\dashrightarrow
\Pbb(V_{\lambda_{i_o}}), f \mapsto f(v_{\lambda_{i_o}})$ and
$r:\Pbb(E)\dashrightarrow \Pbb(V_{\lambda_{i_o}}^*), f \mapsto
v_{\lambda_{i_o}}^*\circ f$. The map $l$ and $r$ are defined on
$\Vcal_{\textrm{End}}$ : they defined respectively $B^-$-equivariant
map from $\Vcal_{\textrm{End}}$ into $\Vcal$, and $B^+$-equivariant
map from $\Vcal_{\textrm{End}}$ into $\Vcal^*$.

The orbits $K_\C\cdot v_{\lambda_{i_o}}\subset
\Pbb(V_{\lambda_{i_o}})$ and $K_\C\cdot v_{\lambda_{i_o}}^*\subset
\Pbb(V_{\lambda_{i_o}}^* )$ are closed and we have
\begin{eqnarray*}
  K_\C\cdot v_{\lambda_{i_o}}\cap \Vcal &=& N^-\cdot v_{\lambda_{i_o}}\simeq N^- \\
  K_\C\cdot v_{\lambda_{i_o}}^*\cap\Vcal^* &=& N^+\cdot v_{\lambda_{i_o}}^*\simeq
  N^+.
\end{eqnarray*}
The rational map $(l,r): \Pbb(E)\dashrightarrow
\Pbb(V_{\lambda_{i_o}})\times \Pbb(V_{\lambda_{i_o}}^*)$ induced
then a map $q : \Vcal_{\textrm{End}}\cap\Xcal_P\to N^-\times N^+$
which is $N^-\times N^+$-equivariant :
$$
q\left( (n^-,n^+)\cdot x\right)= (n^-,n^+)\times q(x)
$$
for $x\in \Vcal_{\textrm{End}}\cap\Xcal_P$, and $n^\pm\in N^\pm$.

We can now finish the arguments. The set $N^- T_\C N^+\subset K_\C$
is dense in $K_\C$, so it is now easy to see that the map
\begin{eqnarray*}
  N^-\times N^+\times \Ycal_P\cap\Vcal_{\textrm{End}} &\longrightarrow &
  \Xcal_P\cap\Vcal_{\textrm{End}} \\
  (n^-,n^+,y) &\longmapsto & (n^-,n^+)\cdot y
\end{eqnarray*}
is a diffeomorphism. We have proved previously that
$\Ycal_P\cap\Vcal_{\textrm{End}}$ is a smooth affine variety, hence
$\Xcal_P$ is smooth near the closed orbit $K_\C\times K_\C\cdot
[v_{\lambda_i}\otimes v_{\lambda_i}^*]\subset
\Xcal\cap\Vcal_{\textrm{End}}$. Lemma \ref{lem:variety-smooth} tells
us then that $\Xcal_P$ is smooth.

\subsection{Hamiltonian actions}

First consider an Hermitian vector space $V$. The Hermitian
structure on $\End(V)$ is $(A,B):= \Tr(AB^*)$, hence the associated
symplectic struture on $\End(V)$ is defined by the relation
$\Omega_{\mathrm{End}}(A,B):=-\textrm{Im}(\Tr(AB^*))$.

Let $U(V)$ be the unitary group. Let $\ugot(V)$ be the Lie algebra
of $U(V)$. We will use the identification
$\epsilon:\ugot(V)\backsimeq \ugot(V)^*$, $X\mapsto \epsilon_X$
where $\epsilon_X(Y)=-\Tr(XY)$. The action $U(V)\times U(V)$ on
$\End(V)$ is $(g,h)\cdot A=gAh^{-1}$. The moment map relative to
this action is
\begin{eqnarray*}
  \End(V)&\longrightarrow& \ugot(V)^*\times\ugot(V)^* \\
  \nonumber  A &\longmapsto&\frac{-1}{2}\left(iAA^*,-iA^*A\right).
\end{eqnarray*}

We consider now the projective space $\Pbb(\End(V))$ equipped with
the Fubini-Study symplectic form $\Omega_{\mathrm{FS}}$. Here the
action of $U(V)\times U(V)$ on $\Pbb(\End(V))$ is hamiltonian with
moment map
\begin{eqnarray*}
   \Pbb(\End(V))&\longrightarrow& \ugot(V)^*\times\ugot(V)^* \\
  \nonumber[A] &\longmapsto&\left(\frac{iAA^*}{\|A\|^2},\frac{-iA^*A}{\|A\|^2}\right).
\end{eqnarray*}

where $\|A\|^2=\Tr(AA^*)$ (see \cite{GIT}[Section 7]). If $\rho:
K\croc U(V)$ is a connected Lie subgroup, we can consider the action
of $K\times K$ on $\Pbb(\End(V))$. Let $\pi_K: \ugot(V)^*\to\kgot^*$
be the projection which is dual to the inclusion $\rho:\kgot\croc
\ugot(V)$. The moment map for the action of $K\times K$ on
$(\Pbb(\End(V)),\Omega_{\mathrm{FS}})$ is then
\begin{eqnarray}
  \Pbb(\End(V))&\longrightarrow & \kgot^*\times\kgot^*  \\
  \nonumber  [A] &\longmapsto &\frac{1}{\|A\|^2}(\pi_K(iAA^*),-\pi_K(iA^*A)).
\end{eqnarray}

We are interested here respectively in
\begin{itemize}
  \item the projective variety $\Xcal_P\subset\Pbb(\End(V))$ with the action of
  $K\times K$,
  \item the projective variety
  $\Ycal_P\subset\Pbb(\End(V))$ with the action of
  $T\times T$,
\end{itemize}
 where $V=\oplus_{i=1}^N V_{\lambda_i}$. The Fubini-Study two-form restrict into
 symplectic forms on $\Xcal_P$
 and $\Ycal_P$. The action of $K\times K$ on $\Xcal_P$ is Hamiltonian
with moment map
 \begin{eqnarray}\label{eq:map-K}
   \Phi_{K\times K}:\Xcal_P &\longrightarrow & \kgot^*\times\kgot^*  \\
  \nonumber  [x] &\longmapsto
  &\frac{1}{\|x\|^2}(\pi_K(ixx^*),-\pi_K(ix^*x)).
\end{eqnarray}

\medskip

Since the diagonal $Z=\{(t,t)| t\in T\}$ acts trivially on $\Ycal_P$
we restrict ourself to the action of $T\times T/Z\simeq T$ on
$\Ycal_P$. Let us compute the moment map $\Phi_{T} : \Ycal_P\to
\tgot^*$ associated to this action. First we have
\begin{equation}\label{eq:phi-T}
    \Phi_{T}([y])=\frac{\pi_T(iy^*y)}{\|y\|^2}=\frac{\pi_T(iyy^*)}{\|y\|^2}
\end{equation}
where $\pi_T: \ugot(V)^*\to\tgot^*$ is the projection which is dual
to $\rho:\tgot\to \ugot(V)$. Since $\rho(X)=i\sum_{j\in
J}\alpha_j(X) v_j\otimes v_j^*$, a small computation shows that for
$B\in \ugot(V)\simeq \ugot(V)^*$ we have $\pi_T(B)=-i\sum_{j\in
J}(Bv_j,v_j)\alpha_j$. Finally for any $[y]\in\Ycal_P$ we get
\begin{equation*}
\Phi_{T}([y])=\sum_{j\in J}\frac{\|yv_j\|^2}{\|y\|^2}\alpha_j.
\end{equation*}

Together with the action on $T$, we have also an action of the Weyl
group $W=N(T)/T$ on $\Ycal_P$ : for $\bar{w}\in W$ we take
\begin{equation}\label{eq:W_action}
\bar{w}\cdot[y]=[\rho(w)\circ y\circ \rho(w)^{-1}],\quad [y]\in
\Ycal_P.
\end{equation}
This action is well defined since the diagonal $Z\subset T\times T$
acts trivially on $\Ycal_P$. The set of weights $\{\alpha_j, j\in
J\}$ is stable under the action of $W$, hence it is an easy fact to
verify that the map $\Phi_T$ is $W$-equivariant.

A dense part of $\Ycal_P$ is formed by the elements $e^Z\cdot[Id]=[\rho(e^{Z})]$.
Take $Z=X+iY\, \in \,\tgot_\C$. We have $\Phi_{T}(e^Z\cdot[Id])=\psi_T(Y)\in
\tgot^*$ with
\begin{equation}\label{eq:phi-rho-T}
    \psi_T(Y)=\frac{1}{\sum_{j\in J} e^{-2\langle\alpha_j,Y\rangle}}
    \sum_{j\in J} e^{-2\langle\alpha_j,Y\rangle}\alpha_j.
\end{equation}
Hence the image of the moment map $\Phi_{T}:\Ycal_P\to \tgot^*$ is
equal to the closure of the image of the map
$\psi_T:\tgot\to\tgot^*$.

\begin{prop}\label{prop:legendre-T}
    The map $\psi_T$ realises a
    diffeomorphism between $\tgot$ and the interior of the polytope
    $P\subset \tgot^*$.
\end{prop}


\textsc{Proof}. Consider the function $F_T:\tgot\to \R$,
$F_T(Y)=\ln\left(\sum_j e^{\langle\alpha_j,Y\rangle}\right)$, and
let $L_{T}:\tgot\to\tgot^*$ be its Legendre transform :
$L_{T}(X)=dF_T|_X$. Note that we have $L_{T}(-2Y)=\psi_T(Y)$.

We see that $F_T$ is strictly convex. So, it is a classical fact
that $L_{T}$ realizes a diffeomorphism of $\tgot$ onto its image, and
 for $\xi\in\tgot^*$ we have
\begin{eqnarray*}\label{eq:image-L}
    \xi\in \textrm{Image}(L_T)&\Leftrightarrow&
    \lim_{Y\to\infty}F_T(Y)-\langle\xi,Y\rangle =\infty\\
    &\Leftrightarrow& \lim_{Y\to\infty}
    \sum_{j\in J} e^{\langle\alpha_j-\xi,Y\rangle}=\infty.
\end{eqnarray*}
In order to conclude we need the following
\begin{lem}\label{lem:lim-F-xi}
    Let $\{\beta_j, j\in J\}$ be a sequence of elements of $\tgot^*$,
    and let $Q$ be its convex hull. We have
$$
\lim_{Y\to\infty}\sum_{j\in J} e^{\langle\beta_j,Y\rangle}=\infty\
\Longleftrightarrow \ 0 \ \in\ \mathrm{Interior}(Q)
$$
\end{lem}

\textsc{Proof}. First we see that $0 \ \notin\ \mathrm{Interior}(Q)$
if and only there exists $v\in\tgot-\{0\}$ such that
$\langle\beta_j,v\rangle\leq 0$ for all $j$ : for such vector $v$,
the map $t\to \sum_{j\in J} e^{t\langle\beta_j,v\rangle}$ is
bounded. Suppose now that $\lim_{Y\to\infty}\sum_{j\in J}
e^{\langle\beta_j,Y\rangle}\neq\infty$. Then there exists a sequence
$(X_k)_k\in\tgot$ such that $\lim_k|X_k|=\infty$ and for all $j$ the
sequence $(\langle\beta_j,X_k\rangle)_k$ remains bounded. If $v$ is
a limit of a subsequence of $(\frac{X_k}{|X_k|})_k$ we have then
$\langle\beta_j,v\rangle\leq 0$ for all $j$. $\Box$

\bigskip

\begin{lem}\label{lem:phi-K-T}
For $[y]\in \Ycal_P$ we have $\Phi_{K\times
K}([y])=(\Phi_{T}([y]),-\Phi_{T}([y]))$.
\end{lem}

\textsc{Proof}. It's sufficient to consider the case
$y=\rho(e^Z)=\sum_{j\in J} e^{i\langle\alpha_j,Z\rangle}v_j\otimes v_j^*$,
for $Z=X+iY\in\tgot_\C$. Then
$yy^*=y^*y=\sum_{j}e^{-2\langle\alpha_j,Y\rangle}v_j\otimes v_j^*=\rho(e^{2iY})$.
So it remains to prove that $\pi_K(iyy^*)=\pi_T(iyy^*)$. We have to check
that $\langle \pi_K(iyy^*), [U,V]\rangle=0$
for $U\in\tgot$ and $V\in\kgot$. We have
\begin{eqnarray*}
  \langle \pi_K(iyy^*), [U,V]\rangle &=& -i\ \Tr\Big(yy^* \rho([U,V])\Big) \\
   &=& -i\ \Tr\Big(\rho(e^{2iY}) [\rho(U),\rho(V)]\Big) \\
   &=& -i\ \Tr\Big([\rho(e^{2iY}),\rho(U)]\rho(V)\Big)= 0.
\end{eqnarray*}

$\Box$

\begin{theo}\label{theo:hamiltonian-X-Y}
We have
\begin{itemize}
      \item $\image(\Phi_{T})=P$,
      \item $\image(\Phi_{K\times K})=
      \{(k_1\cdot\xi,-k_2\cdot\xi)\ |\  \xi\in P \ \mathrm{and} \ k_1,k_2\in K\}$,
      \item $\Ycal_P\subset \Phi_{K\times K}^{-1}(\tgot^*\times \tgot^*)$,
      \item $\Phi_{K\times K}^{-1}(\mathrm{interior}(\Ccal))\subset \Ycal_P$,
      where $\Ccal=C_K\times -C_K$.
\end{itemize}
\end{theo}


\textsc{Proof}. The first point follows from Proposition
\ref{prop:legendre-T}. Since the map $(k_1,t,k_2)\mapsto k_1 t k_2$
from  $K\times T_\C\times K$ into $K_\C$ is onto, we have
\begin{equation}\label{eq:X-Y}
    \Xcal_P= (K\times K)\cdot \Ycal_P.
\end{equation}
So if $[x]\in \Xcal_P$, there exist $[y]\in\Ycal$ and $k_1,k_2\in K$ such
that $[x]=(k_1,k_2)\cdot [y]$, hence
\begin{eqnarray}\label{eq:phi-K-T}
\nonumber \Phi_{K\times K}([x])&=&(k_1,k_2)\cdot \Phi_{K\times K}([y])\\
&=& \left(k_1\cdot\Phi_{T}([y]),-k_2\cdot\Phi_{T}([y])\right)
\end{eqnarray}
The second point is proved. The third point follows also from the
identity (\ref{eq:phi-K-T}) when $k_1=k_2=e$. Consider now
$[x]=(k_1,k_2)\cdot [y]$ such that $\Phi_{K\times K}([x])$ belongs
to the interior of the cone $C_K\times -C_K$. Then
$k_1\cdot\Phi_{T}([y])$ and $k_2\cdot\Phi_{T}([y])$ are regular
points of $C_K$. This implies that $k_1,k_2\in N(T)$ and
$k_2k_1^{-1}\in T$. So
\begin{eqnarray*}
[x]&=&(k_1,k_2)\cdot [y]\\
   &=&(e,k_2k_1^{-1})\cdot \Big((k_1,k_1)\cdot [y]\Big)\in \Ycal_P
\end{eqnarray*}
since $\Ycal_P$ is stable under the actions of $T\times T$ and $W$.
$\Box$

\begin{coro}\label{coro:vertex-P}
If $[v_{\lambda_i}\otimes v_{\lambda_i}^*]\in \Xcal_P$ then
$\lambda_i$ is a vertex of the polytope $P$.
\end{coro}

\textsc{Proof}. Let $x=v_{\lambda_i}\otimes v_{\lambda_i}^*$, and suppose that
$[x]$ belongs to $\Xcal_P$. In order to show that $[x]\in \Ycal_P$, we compute
$\Phi_{K\times K}([x])$. We see that $xx^*=x^*x=x$ and $\|x\|=1$ so
$\Phi_{K\times K}([x])=(\pi_K(ix),-\pi_K(ix))$. For $X\in\kgot$ we have
\begin{eqnarray*}
\langle\pi_K(ix),X\rangle &=& -i\ \Tr\Big(v_{\lambda_i}\otimes v_{\lambda_i}^* \rho(X)\Big) \\
   &=& -i\ (\rho(X)v_{\lambda_i},v_{\lambda_i})\\
   &=&  \langle\lambda_i,X\rangle.
\end{eqnarray*}
We have then $\Phi_{K\times K}([x])=(\lambda_i,-\lambda_i)$ with
$\lambda_i$ beiing a regular point of $C_K$ : hence $[x]\in
\Ycal_P$. Now we can conclude with the help of Lemma
\ref{lem:Y-closed-orbit}. Since $[v_{\lambda_i}\otimes
v_{\lambda_i}^*]$  belongs to $\Ycal_P$, the weight $\lambda_i$ is a
vertex of the polytope $P$. $\Box$

\begin{rem}\label{rem:smoothness}
In this section, Theorem \ref{theo:hamiltonian-X-Y} was obtain
without using the fact that the varieties $\Xcal_P$ and $\Ycal_P$
are smooth. Hence Corollary \ref{coro:vertex-P} can be used to prove
the smoothness of $\Xcal_P$.
\end{rem}

\subsection{Symplectic cutting}\label{subsec:cutting}

Let $(M,\Omega_M,\Phi_M)$ be an Hamiltonian $K$-manifold. At this
stage the moment map $\Phi_M$ is \emph{not assumed to be proper}. We
consider also the Hamiltonian $K\times K$-manifold $\Xcal_P$
associated to a $K$-adapted polytope $P$.

The purpose of this section is to define a {\em symplectic cutting}
of $M$ which uses $\Xcal_P$. The notion of symplectic cutting was
introduced by Lerman in \cite{Lerman-cut} in the case of a torus
action. Later Woodward \cite{Woodward96} extends this procedure to
the case of a non-abelian group action (see also
\cite{Meinrenken98}). The method of symplectic cutting that we
define in this section is different from the one of Woodward.

\medskip

We have two actions of $K$ on $\Xcal_P$ : the action from the left
(resp. right) , denoted $\cdot_l$ (resp. $\cdot_r$), with moment map
$\Phi_l:\Xcal_P\to \kgot^*$ (resp. $\Phi_r$). We consider now the
product of $M\times \Xcal_P$ with
\begin{itemize}
  \item the action $k\cdot_1 (m,x)=(k\cdot m,k\cdot_r x)$ : the
  corresponding moment map is $\Phi_1(m,x)=\Phi_M(m)+\Phi_r(x)$,
  \item the action  $k\cdot_2 (m,x)=(m,k\cdot_l x)$ : the
  corresponding moment map is $\Phi_2(m,x)=\Phi_l(x)$.
\end{itemize}

\begin{defi}\label{def:M-P}
    We denote $M_P$ the symplectic reduction at $0$ of $M\times \Xcal_P$
    for the action $\cdot_1$ : $M_P:=(\Phi_1)^{-1}(0)/K$.
\end{defi}

Note that $M_P$ is compact when $\Phi_M$ is proper. The action
$\cdot_2$ on $M\times \Xcal_P$ induces an action of $K$ on $M_P$.
The moment map $\Phi_2$ induces an equivariant map $\Phi_{M_P}:
M_P\to \kgot^*$. Let $\Zcal\subset (\Phi_1)^{-1}(0)$ be the set of
points where $(K,\cdot_1)$ as a trivial stabilizer.

\begin{defi}\label{def:smooth}
    We denote $M_P'$ the quotient $\Zcal/K\subset M_P$.
\end{defi}

$M_P'$ is an open subset of smooth points of $M_P$ which is
invariant under the $K$-action. The symplectic structure of $M\times
\Xcal_P$ induces a canonical symplectic structure on $M_P'$ that we
denote $\Omega_{M_P'}$. The action of $K$ on $(M_P', \Omega_{M_P'})$
is Hamiltonian with moment map equal to the restriction of
$\Phi_{M_P}: M_P\to \kgot^*$ to $M_P'$.

We start with the easy

\begin{lem}\label{lem:image-phi-M-P}
    The image of $\Phi_{M_P}: M_P\to \kgot^*$ is equal to
    the intersection of the image of $\Phi_M: M\to\kgot^*$ with $K\cdot P$.
\end{lem}

\medskip

Let $\Ucal_P=K\cdot\mathrm{Interior}(P)\subset K\cdot P$. We will
show now that the open and dense subset $(\Phi_{M_P})^{-1}(\Ucal_P)$
of $M_P$ belongs to $M_P'$. Afterwards we will prove that
$\Phi_{M_P}^{-1}(\Ucal_P)$ is quasi-symplectomorphic to the open
subset $\Phi_{M}^{-1}(\Ucal_P)$ of $M$.

\medskip

 We consider the open and dense subset of
$\Xcal_P$ which is equal to the open orbit $\bar{\rho}(K_\C)$. From
Lemma \ref{lem:properties-rho}, we know that
\begin{eqnarray}\label{eq:Theta}
    \Theta : K\times \kgot & \longrightarrow &  \bar{\rho}(K_\C)\\
\nonumber (k,X)            & \longmapsto &   [\rho(ke^{iX})]
\end{eqnarray}
is a diffeomorphism. Through $\Theta$, the action of $K\times K$ on
$K\times \kgot$ is $k\cdot_l(a,X)=(ka,X)$ for the action "from the
left" and $k\cdot_r(a,X)=(ak^{-1},k\cdot X)$ for the action "from
the right".

We consider now the map $\psi_K :\kgot \to \kgot^*$ defined by
$\psi_K(X)=\Phi_l([\rho(e^{iX})])$. In other words,
$$
\psi_K(X)=\frac{\pi_K(i\rho(e^{i2X}))}{\Tr(\rho(e^{i2X}))}.
$$

Consider the function $F_K:\kgot\to \R$,
$F_K(X)=\ln(\Tr(\rho(e^{-iX}))$. Let $L_{K}:\kgot\to\kgot^*$ be its
Legendre transform.

\begin{prop}\label{prop:legendre-K}
\begin{itemize}
  \item We have $\psi_K(X)=L_K(-2X)$, for $X\in\kgot$,
  \item The function $F_K$ is strictly convex,
  \item The map $\psi_K$ realizes an equivariant diffeomorphism between $\kgot$ and $\Ucal_P$.
  \item The image of $\Phi_l: \Xcal_P\to\kgot^*$ is equal to the closure
    of $\Ucal_P$,
  \item   $\Phi_l^{-1}(\Ucal_P)=\bar{\rho}(K_\C)$.
\end{itemize}
\end{prop}

\textsc{Proof}. For $X,Y\in\kgot$ we consider the function $\tau(s)=F_K(X+sY)$.
Since $F_K$ is $K$-invariant
we can restrict our computation to $X\in\tgot$. We will use the decomposition of
$Y\in\kgot$ relatively to the $T$-weights on $\kgot_\C$ : $Y=\sum_\alpha Y_\alpha$
where $\ad(Z)Y_\alpha=i \alpha(Z) Y_\alpha$ for any $Z\in \tgot$, and $Y_0\in\tgot$.
We have
\begin{eqnarray*}
\tau'(s)&=&\frac{-i}{\Tr(\rho(e^{-iX_s}))}
\Tr\left(\rho(e^{-iX_s}) \rho\left(\frac{e^{i\ad(X_s)}-1}{i\ad(X_s)}Y\right)\right)\\
&=&\frac{-i}{\Tr(\rho(e^{-iX_s}))}
\Tr\left(\rho(e^{-iX_s}) \rho(Y)\right)\\
   &=& \frac{1}{\Tr(\rho(e^{-iX_s}))}\langle\pi_K(i\rho(e^{-iX_s})),Y\rangle
\end{eqnarray*}
where $X_s =X+sY$. Since by definition $\tau'(0)=\langle L_K(X),Y\rangle$,
the first point is proved. For the second derivative we have
\begin{eqnarray*}
  \tau''(0) &=& -\left(\frac{\Tr(\rho(e^{-iX}) \rho(iY))}{\Tr(\rho(e^{-iX}))}\right)^2
+ \frac{\Tr\left(\rho(e^{-iX}) \rho(\frac{e^{i\ad(X)}-1}{i\ad(X)}iY)\rho(iY)\right)}
{\Tr(\rho(e^{-iX}))}\\
   &=& R_1 + R_2
\end{eqnarray*}
where
\begin{eqnarray*}
R_1&=& \frac{\Tr\left(\rho(e^{-iX}) \rho(iY_0)\rho(iY_0)\right)}
{\Tr(\rho(e^{-iX}))}-\left(\frac{\Tr(\rho(e^{-iX}) \rho(iY_0))}{\Tr(\rho(e^{-iX}))}\right)^2\\
&=&\frac{\sum_j e^{-\langle\alpha_j,X\rangle}\langle\alpha_j,Y_0\rangle^2}
{\sum_j e^{-\langle\alpha_j,X\rangle}}-
\left(\frac{\sum_j e^{-\langle\alpha_j,X\rangle}\langle\alpha_j,Y_0\rangle}
{\sum_j e^{-\langle\alpha_j,X\rangle}}\right)^2
\end{eqnarray*}
and
\begin{eqnarray*}
R_2&=& \frac{1}{\Tr(\rho(e^{-iX}))}\sum_{\alpha\neq 0,\beta\neq 0}
\frac{e^{-\langle\alpha,X\rangle}-1}{-\langle\alpha,X\rangle}
\Tr\left(\rho(e^{-iX}) \rho(iY_\alpha)\rho(iY_\beta)\right)\\
&=& \frac{1}{\Tr(\rho(e^{-iX}))}\sum_{\alpha\neq 0,j}
\frac{e^{-\langle\alpha,X\rangle}-1}{-\langle\alpha,X\rangle}e^{-\langle\alpha_j,X\rangle}
\|\rho(Y_\alpha)v_j\|^2.
\end{eqnarray*}
It is now easy to see that $R_1$ and $R_2$ are positive and that $R_1+R_2>0$ if $Y\neq 0$.
We have proved that $F_K$ is strictly convex, So, its Legendre transform
$L_{K}$ realizes a diffeomorphism of $\kgot$ onto its image. Using the first point we
know that $\psi_K$ realizes a diffeomorphism of $\kgot$ onto its image.
The map $\psi_K$ is equivariant and coincides with $\psi_T$ on
$\tgot$. We have proved in Proposition \ref{prop:legendre-T} that the image of
$\psi_T$ is equal to the interior of $P$, hence the image of $\psi_K$ is $\Ucal_P$.

For the last two points we first remark that
\begin{equation}\label{eq-phipsi}
    \Phi_l([\rho(k e^{iX})])=k\cdot \psi_K(X)
\end{equation}
hence the image of $\Phi_l$ is the closure of $\Ucal_P$. If we use
the fact that $\psi_K$ is a diffeomorphism from $\kgot$ onto
$\Ucal_P$,(\ref{eq-phipsi}) shows that $\Phi_l^{-1}(K\cdot\xi)\cap
\bar{\rho}(K_\C)$ is a non empty and closed subset of
$\Phi_l^{-1}(K\cdot\xi)$ for any $\xi\in\Ucal_P$ (in fact it is a
$K\times K$-orbit). On the other hand $\Phi_l^{-1}(K\cdot\xi)\cap
(\Xcal_P \setminus \bar{\rho}(K_\C))$ is also a closed subset of
$\Phi_l^{-1}(K\cdot\xi)$ since $\bar{\rho}(K_\C))$ is open in
$\Xcal_P$. Since $\Phi_l^{-1}(K\cdot\xi)$ is connected the second
subset is empty : in other words $\Phi_l^{-1}(K\cdot\xi)\subset
\bar{\rho}(K_\C)$. $\Box$

\bigskip

We introduce now the equivariant diffeomorphism
\begin{eqnarray}\label{eq:upsilon}
    \Upsilon  : K\times \Ucal_P & \longrightarrow &  \bar{\rho}(K_\C)\\
\nonumber (k,\xi)   & \longmapsto &\Theta(k,\psi^{-1}_K(\xi)).
\end{eqnarray}
We now look at $K\times \Ucal_P$ equipped with the symplectic
structure $\Upsilon^*(\Omega_{\Xcal_P})$, and the Hamiltonian action
of $K\times K$ : the moment maps satisfy
\begin{equation}\label{eq:moment-map-upsilon}
    \Upsilon^*(\Phi_l)(k,\xi)= k\cdot \xi\quad and \quad
\Upsilon^*(\Phi_r)(k,\xi)=-\xi.
\end{equation}

\begin{prop}\label{prop:structure-K-U}
    We have
$$
\Upsilon^*(\Omega_{\Xcal_P})=d\lambda + d\eta
$$
where $\lambda$ is the Liouville $1$-form on $K\times\kgot^*\simeq
\T^*K$ and $\eta$ is an invariant $1$-form on
$\Ucal_P\subset\kgot^*$ which is killed by the vectors tangent to
the $K$-orbits.
\end{prop}

\textsc{Proof}. Let $E_1,\dots,E_r$ be a basis of $\kgot$, with dual
basis $\xi^1,\dots,\xi^r$. Let $\omega^i$ the $1$-form on $K$,
invariant by left translation and  equal to $\xi^i$ at the identity.
The Liouville $1$-form is $\lambda=-\sum_i \omega^i\otimes E_i$. For
$X\in\kgot$ we denote $X_l(k,\xi)=\frac{d}{dt}|_0
e^{-tX}\cdot_l(k,\xi)$ and $X_r(k,\xi)=\frac{d}{dt}|_0
e^{-tX}\cdot_r(k,\xi)$ the vectors fields generated by the action of
$K\times K$. Since $\iota(X_l)d\lambda=-d\langle\Phi_l,X\rangle$ and
$\iota(X_r)d\lambda=-d\langle\Phi_r,X\rangle$, the closed invariant
$2$-form $\beta=\Upsilon^*(\Omega_{\Xcal_P})-d\lambda$ is $K\times
K$ invariant and is killed by the vectors tangent to the orbits :
$(*) \,$ $\iota(X_l)\beta=\iota(X_r)\beta=0$ for all $X\in\kgot$. We
have $\beta=\beta_2+\beta_1+\beta_0$ where
$\beta_2=\sum_{i,j}a_{ij}(\xi)\omega^i\wedge\omega^j$,
$\beta_1=\sum_{i,j}b_{ij}(\xi)\omega^i\wedge dE_j$, and $\beta_0$ is
an invariant $2$-form on $\Ucal_P$. The equalities $(*)$ gives
$\iota(X_l)\beta_2=\iota(X_l)\beta_1=0$ which implies that
$\beta_2=\beta_1=0$. So $\beta=\beta_0$ is a closed invariant
$2$-form on $\Ucal_P$ which is killed by the vectors tangent to the
$K$-orbits. Since $\Ucal_P$ admits a retraction to $\{0\}$,
$\beta=d\eta$ where $\eta$ is an invariant $1$-form on $\Ucal_P$
which is killed by the vectors tangent to the $K$-orbits. $\Box$

\bigskip

If $(m,x)\in M\times \Xcal_P$ belongs to $\Phi_{1}^{-1}(0)$, we
denote $[m,x]$ the corresponding element in $M_P$. By definition we
have $\Phi_{M_P}([m,x])=\Phi_l(x)$ for $[m,x]\in M_P$, hence the
image of $\Phi_{M_P}$ is included in the closure of $\Ucal_P$. We
see also that $[m,x]\in \Phi_{M_P}^{-1}(\Ucal_P)$ if and only if
$x\in \Phi_l^{-1}(\Ucal_P)=\bar{\rho}(K_\C)$. Since  $(K,\cdot_r)$
acts freely on $\bar{\rho}(K_\C)$, we see that $(K,\cdot_1)$ acts
freely on $\Phi_{M_P}^{-1}(\Ucal_P)$ : the open and dense set
$\Phi_{M_P}^{-1}(\Ucal_P)\subset M_P$ is then contained in  $M_P'$.

Now, we can state our main result which compares the open invariant
subsets $\Phi_M^{-1}(\Ucal_P)\subset M$ and
$\Phi_{M_P}^{-1}(\Ucal_P)\subset M_P$ equipped respectively with the
symplectic structures $\Omega_M$ and $\Omega_{M_P'}$.

\begin{theo}\label{theo:M-P-smooth}
$\Phi_{M_P}^{-1}(\Ucal_P)$ is an open and dense subset of smooth
points in $M_P$. There exist an equivariant diffeomorphism
$\Psi:\Phi_M^{-1}(\Ucal_P)\to\Phi^{-1}_{M_P}(\Ucal_P)$ such that
$$
\Psi^*(\Omega_{M_P'})= \Omega_M + d\Phi_M^*\eta.
$$
Here $\eta$ is an invariant $1$-form on $\Ucal_P$ which is killed by
the vectors tangent to the $K$-orbits. Moreover the path $\Omega^t=
\Omega_M + t d\Phi_M^*\eta, $ defines an homotopy of symplectic
$2$-forms between $\Omega_M$ and $\Psi^*(\Omega_{M_P'})$.
\end{theo}

\begin{rem}\label{rem-quasi-symplectomorphism}
    The map $\Psi$ will be call a \emph{quasi}-symplectomorphism.
\end{rem}

\textsc{Proof}.  Consider the immersion
\begin{eqnarray*}
  \psi:  \Phi_M^{-1}(\Ucal_P)&\longrightarrow & M\times \Xcal_P \\
   m &\longmapsto & (m,\Upsilon(e,\Phi_M(m))).
\end{eqnarray*}
We have $\Phi_1(\psi(m))=\Phi_M(m)+
\Upsilon^*\Phi_r(e,\Phi_M(m))=0$, and $\Phi_2(\psi(m))=$ \break
$\Upsilon^*\Phi_l(e,\Phi_M(m))=\Phi_M(m)\in\Ucal_P$ (see
(\ref{eq:moment-map-upsilon})). Hence for all
$m\in\Phi_M^{-1}(\Ucal_P)$, we have $\psi(m)\in\Phi_{1}^{-1}(0)$,
and its class $[\psi(m)]\in M_P$ belongs to
$\Phi_{M_P}^{-1}(\Ucal_P)$.

We denote $\Psi:\Phi_M^{-1}(\Ucal_P)\to\Phi^{-1}_{M_P}(\Ucal_P)$ the
map $m\mapsto [\psi(m)]$. Let us show that it defines a
diffeomorphism. If $\Psi(m)=\Psi(m')$, there exists $k\in K$ such
that
\begin{eqnarray*}
(m,\Upsilon(e,\Phi_M(m)))&=&k\cdot_1(m',\Upsilon(e,\Phi_M(m')))\\
&=& (k\cdot m',k\cdot_r\Upsilon(e,\Phi_M(m')))\\
&=& (k\cdot m',\Upsilon(k^{-1},k\cdot\Phi_M(m'))).
\end{eqnarray*}
Since $\Upsilon$ is a diffeomorphism, we must have $k=e$ and $m=m'$
: the map $\Psi$ is one to one. Consider now $(m,x)\in
\Phi_{1}^{-1}(0)$ such that $\Phi_{M_P}([m,x])=\Phi_l(x)\in \Ucal_P$
: then $x\in\Phi_l^{-1}(\Ucal_P)=\bar{\rho}(K_\C)=\image(\Upsilon)$.
We have $x=\Upsilon(k,\xi)$ where $\xi=-\Phi_r(x)=\Phi_M(m)$.
Finally
\begin{eqnarray*}
(m,x)&=&(m,\Upsilon(k,\Phi_M(m)))\\
&=& k^{-1}\cdot_1(k\cdot m,\Upsilon(e,k\cdot\Phi_M(m)))\\
 &=& k^{-1}\cdot_1\psi(k\cdot m).
\end{eqnarray*}
We have proved that $\Psi$ is onto.

In order to show that $\Psi$ is a submersion we must show that for
$m\in\Phi_M^{-1}(\Ucal_P)$
$$
\image(\T_m\psi)\oplus \T_{\psi(m)}(K\cdot_1\psi(m))=
\T_{\psi(m)}\Phi_{1}^{-1}(0).
$$
Here $\T_m\psi: \T_m M\to \T_{\psi(m)}(M\times \Xcal_P)$ is the
tangent map, and $\T_{\psi(m)}(K\cdot_1\psi(m))$ denotes the tangent
space at $\psi(m)$ of the $(K,\cdot_1)$-orbit. We have
$\dim(\image(\T_m\psi))+ \dim(\T_{\psi(m)}(K\cdot_1\psi(m)))=
\dim(\T_{\psi(m)}\Phi_{1}^{-1}(0))$ so it is sufficient to prove
that
$$
\image(\T_m\psi)\cap\T_{\psi(m)}(K\cdot_1\psi(m))=\{0\}.
$$
Consider $(v,w)\in
\image(\T_m\psi)\cap\T_{\psi(m)}(K\cdot_1\psi(m))$. There exists
$X\in\kgot$ such $(v,w)=\frac{d}{dt}|_0 e^{tX}\cdot_1 \psi(m)$ :
$$
v=\frac{d}{dt}_{|_0} e^{tX}\cdot m \quad \mathrm{and}\quad
w=\frac{d}{dt}_{|_0} e^{tX}\cdot_r\Upsilon(e,\Phi_M(m))
$$
In the other hand since $(v,w)\in \image(\T_m\psi)$, we have
$$
w=\frac{d}{dt}_{|_0} \Upsilon(e,\Phi_M(e^{tX}\cdot m))
$$
Since
$e^{tX}\cdot_r\Upsilon(e,\Phi_M(m))=\Upsilon(e^{-tX},\Phi_M(e^{tX}\cdot
m))$ we obtain that
$$
\frac{d}{dt}_{|_0}\Upsilon(e^{-tX},\Phi_M(e^{tX}\cdot m))=
\frac{d}{dt}_{|_0} \Upsilon(e,\Phi_M(e^{tX}\cdot m))
$$
or in other words $\frac{d}{dt}_{|_0}\Upsilon(e^{-tX},\Phi_M(m))=0$.
Since $\Upsilon$ is a diffeomorphism we have $X=0$, and then
$(v,w)=0$.

We can now compute the pull-back by $\Psi$ of the symplectic form
$\Omega_{M_P'}$. We have
\begin{eqnarray*}
  \Psi^*(\Omega_{M_P'})&=& \psi^*(\Omega_M + \Omega_{X_P})\\
   &=& \Omega_M + \Phi_M^*\Upsilon^*(\Omega_{\Xcal_P})\\
   &=& \Omega_M + d\Phi_M^*\eta.
\end{eqnarray*}

It remains to prove that for every $t\in [0,1]$, the $2$-form
$\Omega^t= \Omega_M + t d\Phi_M^*\eta$ is non-degenerate. Take
$t\neq 0$, $m\in \Phi_M^{-1}(\Ucal_P)$ and suppose that the
contraction of $\Omega^t|_m$ by $v\in \T_m M$ is equal to $0$. For
every $X\in \kgot$ we have
\begin{eqnarray*}
  0 &=& \Omega^t(X_M(m),v) \\
   &=& -\iota(v) d\langle\Phi_M,X\rangle_{|_m}+ t\iota(v)
   \iota(X_M)d\Phi_M^*\eta_{|_m}\\
   &=&-\iota(v) d\langle\Phi_M,X\rangle_{|_m}
\end{eqnarray*}
since $\iota(X_M)d\Phi_M^*\eta=
d\Phi_M^*(\iota(X_{\kgot^*})\eta)=0$. Thus we have
$\T_m\Phi_M(v)=0$, and then $\iota(v)d\Phi_M^*\eta=0$. Finally we
have that $0=\iota(v)\Omega^t|_m=\iota(v)\Omega_M|_m$. But
$\Omega_M$ is non-degenerate, so $v=0$. $\Box$.

\subsection{Formal quantization : second definition}

We suppose here that the Hamiltonian $K$-manifold
$(M,\Omega_M,\Phi_M)$ is \emph{proper} and admits a Kostant-Souriau
line bundle $L$. Now we consider the complex $K\times K$-submanifold
$\Xcal_P$ of $\Pbb(E)$. Since $\Ocal(-1)$ is a $K\times
K$-equivariant Kostant-Souriau line bundle on the projective space
$\Pbb(E)$ the restriction
\begin{equation}\label{eq:L-P}
    L_P=\Ocal(-1)|_{\Xcal_P}
\end{equation}
is a Kostant-Souriau line bundle on $\Xcal_P$. Hence $L\boxtimes
L_P$ is a Kostant line bundle on the product $M\times \Xcal_P$. In
section \ref{subsec:quant-singular} we have have defined the
quantization $\Qcal_K(M_P)$ of the (singular) reduced space
$M_P:=(M\times \Xcal_P)/\!\!/_0 (K,\cdot_1)$.

\medskip

\textsc{Notation}: $O_K(r)$ will be any element $\sum_{\mu\in
\what{K}}m_\mu V_\mu^K$ of $\Rfor(K)$ where $m_\mu=0$ if $\|\mu\|<
r$. The limit $\lim_{r\to+\infty}O_K(r)=0$ defines the notion of
convergence in $\Rfor(K)$.

\medskip

\begin{prop}\label{prop:Q-M-P}
    Let $\varepsilon_P>0$ be the radius of the biggest ball center
    at $0\in\tgot^*$ which is contains in  the polytope $P$.
    We have
    \begin{equation}\label{eq-Q-MP}
    \Qcal_K(M_P)= \sum_{\|\mu\|<\varepsilon_P}\Qcal(M_\mu)
    V_\mu^K \ +\ O_K(\varepsilon_P).
    \end{equation}
\end{prop}

\textsc{Proof}. Theorem \ref{theo:Q-R-singular} - ``Quantization
commutes with reduction in the singular setting'' - tells us that $
  \Qcal_K(M_P)= \sum_{\mu\in \what{K}}\Qcal((M_P)_\mu)
    V_\mu^K$ where $(M_P)_\mu$ is the symplectic reduction
    $$
    (M_P\times \overline{K_2\cdot\mu})/\!\!/_{0}K_2\cong(M\times \Xcal_P\times
\overline{K_2\cdot\mu})/\!\!/_{(0,0)}K_2\times
    K_1.
    $$
Recall what the $K_1, K_2$-action are:
$k\cdot_1(m,x,\xi)=(km,k\cdot_r x,\xi)$ and $k\cdot_2(m,x,\xi)=
(m,k\cdot_l x,k\xi)$ for $(m,x,\xi)\in M\times \Xcal_P\times
\overline{K_2\cdot\mu}$ and $k\in K$.

Since the image of $\Phi_{M_P}$ is equal to the intersection of
$K\cdot P=\overline{\Ucal_P}$ with the image of $\Phi_M$, we have
\begin{equation}\label{eq:Q-M-P-mu}
    \Qcal((M_P)_\mu)=0\quad \mathrm{if}\quad \mu\notin P\cap
    \mathrm{Image}(\Phi_M).
\end{equation}
We will now exploit Theorem \ref{theo:M-P-smooth} to show that
$\Qcal((M_P)_\mu)=\Qcal(M_\mu)$ if $\mu$ belongs to the interior of
$P$.

There exists a quasi-symplectomorphism $\Psi$ between the open
subset $\Phi^{-1}_{M}(\Ucal_P)$ of $M$ and the open and dense subset
$\Phi^{-1}_{M_P}(\Ucal_P)$ of $M_P$. Moreover one can see easily
that the restriction of the Kostant line bundle $L_P\to \Xcal_P$ to
the open subset $\bar{\rho}(K_\C)$ is \emph{trivial}. If $L_{M_P}$
is the Kostant line bundle on $M_P$ induced by $L\boxtimes L_P$, we
have that the pull-back of the restriction
$L_{M_P}|_{\Phi^{-1}_{M}(\Ucal_P)}$ by $\Psi$ is equivariantly
diffeomorphic to the restriction of $L$ to $\Phi^{-1}_{M}(\Ucal_P)$.

Take now $\mu\in \what{K}$ that belongs to the interior of the
polytope $P$.  The element $\Qcal((M_P)_\mu)\in \Z$ is given by the
index of a transversally elliptic symbol defined in a (small)
neighborhood of $\Phi^{-1}_{M_P}(\mu)\subset M_P$. This symbol is
defined through two auxiliary data: the Kostant line bundle
$L_{M_P}$ and a compatible almost complex structure $J$ which
defined in a neighborhood of $\Phi^{-1}_{M_P}(\mu)$. If we pull back
everything by $\Psi$, we get a transversally elliptic symbol living
in a (small) neighborhood of $\Phi^{-1}_{M}(\mu)\subset M$ which is
defined by the Kostant line bundle $L$ and an almost complex
structure $J_1$ \emph{compatible with the symplectic structure}
$\Omega_1:=\Omega_M + d\Phi_M^*\eta$. But since $\Omega_t=\Omega_M +
td\Phi_M^*\eta$ defines an homotopy of symplectic structures, any
almost complex structure compatible with $\Omega_M$ is homotopic to
$J_1$. We have then shown that $\Qcal(M_\mu)=\Qcal((M_P)_\mu)$ for
any $\mu$ belonging to the interior of $P$. So we have
$$
\Qcal_K(M_P)= \sum_{\mu\in \mathrm{Interior}(P)}\Qcal(M_\mu)
    V_\mu^K + \sum_{\nu\in \partial P}\Qcal((M_P)_\nu)
    V_\nu^K.
$$
Since for $\nu \in \partial P$ we have $\|\nu\|\geq \varepsilon_P$,
the last equality proves (\ref{eq-Q-MP}). $\Box$

\medskip

We work now with the dilated polytope $nP$, for any integer $n\geq
1$. The polytope $nP$ is still $K$-adapted, so one can consider the
reduced spaced $M_{nP}$ and Proposition \ref{prop:Q-M-P} gives that
\begin{equation}\label{eq-Q-MkP}
    \Qcal_K(M_{nP})= \sum_{\|\mu\|<n\varepsilon_P}\Qcal(M_\mu)
    V_\mu^K \ +\ O_K(n\varepsilon_P).
\end{equation}
for any integer $n\geq 1$. We can summarize the result of this
section in the following

\begin{prop}\label{prop:formal-second-definition}
Let $(M,\Omega_M)$ be a pre-quantized Hamitonian $K$-manifold, with
a proper moment map $\Phi_M$.

$\bullet$ For any integer $n\geq 1$, the (singular) compact
Hamiltonian manifold $M_{nP}$ contains as an open and dense subset,
the open subset  $\Phi^{-1}_M(n\Ucal_P)$ of $M$.

$\bullet$ We have $\qfor_K(M)=\lim\limits_{n\to
\infty}\Qcal_K(M_{nP})$.
\end{prop}

\section{Functorial properties : Proof of Theorem \ref{theo:intro}}\label{sec:proof}

This section is devoted to the proof of Theorem \ref{theo:intro}. We
will used in a crucial way the characterisation of $\qfor_K$ given
in Proposition \ref{prop:formal-second-definition}.

\medskip

Let $H\subset K$ be a connected Lie subgroup. Here we consider a
pre-quantized Hamiltonian $K$-manifold $M$ which is \emph{proper} as
an Hamiltonian $H$-manifold. We want to compare $\qfor_K(M)$ and
$\qfor_H(M)$. For $\mu\in \what{K}$ and $\nu\in \what{H}$ we denote
$N^\mu_{\nu}$ the multiplicity of $V_\nu^H$ in the restriction
$V_\mu^K|_H$. We have seen in the introduction that $N^\mu_{\nu}
\Qcal\left(M_{\mu,K}\right)\neq 0$ only for the $\mu$ belonging to
\emph{finite} subset $ \what{K}\cap \Phi_K\left(
K\cdot\Phi^{-1}_H(\nu)\right)$. Then $\qfor_K(M)$ is $H$-admissible
and we have the following equality in $\Rfor(H)$ :
\begin{equation}\label{eq:restriction-formal}
\qfor_K(M)|_H=\sum_{\nu\in \what{H}} m_\nu V_\nu^H
\end{equation}
with $m_\nu=\sum_\mu N^\mu_{\nu} \Qcal\left(M_{\mu,K}\right)$. We
will now prove that $\qfor_K(M)|_H=\qfor_H(M)$.

\medskip

\begin{lem}\label{lem:K-H}
The restriction $\qfor_K(M)|_H$ is equal to
$\lim\limits_{n\to\infty}\Qcal_K(M_{nP})|_H$.
\end{lem}

\textsc{Proof}. Let us denote $P^o$ and $\partial P$ respectively
the interior and the border of the $K$-adapted polytope $P$. We
write
$$
\qfor_K(M)=\sum_{\mu\in n P^o}\Qcal(M_{\mu,K}) V_\mu^K
+\sum_{\mu\notin nP^o}\Qcal(M_{\mu,K}) V_\mu^K.
$$
On the other side
$$
\Qcal_K(M_{nP})=\sum_{\mu\in n P^o}\Qcal(M_{\mu,K}) V_\mu^K +
\sum_{\mu\in n\partial P}\Qcal((M_{nP})_{\mu,K}) V_\mu^K.
$$
So the difference $D(n)=\qfor_K(M)-\Qcal_K(M_{nP})$ is equal to
$$
D(n)=-\sum_{\mu'\in n\partial P}\Qcal((M_{nP})_{\mu',K}) V_\mu^K
+\sum_{\mu\notin nP^o}\Qcal(M_{\mu,K}) V_\mu^K.
$$
We show now that the restriction $D(n)|_H$ tends to $0$ in
$\Rfor(H)$ as $n$ goes to infinity. For this purpose, we will prove
that for any $c>0$ there exist $n_c\in\N$ such that $D(n)|_H=O_H(c)$
for any $n\geq n_c$.

For $c>0$ we consider the compact subset of $\kgot^*$ defined by
\begin{equation}\label{eq:C-c}
    \Kcal_c=\Phi_K\left(K\cdot\Phi^{-1}_H(\xi\in\hgot^*, \|\xi\|\leq
c)\right).
\end{equation}

Let $n_c\in\N$ such that $\Kcal_c$ is included in $K\cdot (n_c P^o)$
: hence $\Kcal_c\subset K\cdot (n P^o)$ for any $n\geq n_c$. We know
that for $\mu\in\what{K}$, we have $N^\mu_{\nu}
\Qcal\left(M_{\mu,K}\right)\neq 0$ only for $\mu\in\Phi_K\left(
K\cdot\Phi^{-1}_H(\nu)\right)$, and for $\mu'\in\what{K}$, we have
$N^{\mu'}_{\nu} \Qcal\left((M_{nP})_{\mu',K}\right)\neq 0$ only for
$\mu'\in nP\cap\Phi_K\left( K\cdot\Phi^{-1}_H(\nu)\right)$.

Then if $n\geq n_c$, we have
$$
N^\mu_{\nu}
\Qcal\left(M_{\mu,K}\right)=N^{\mu'}_{\nu}
\Qcal\left((M_{nP})_{\mu',K}\right)=0
$$
for any $\nu\in \what{H}\cap\{\xi\in\hgot^*, \|\xi\|\leq c\}$,
$\mu\notin nP^o$ and $\mu'\in n\partial P$. It means that
$D(n)|_H=O_H(c)$ for any $n\geq n_c$. $\Box$

\medskip

Since $\Qcal_K(M_{nP})|_H=\Qcal_H(M_{nP})$, we are no led to the

\begin{lem}
The limit $\lim\limits_{n\to\infty}\Qcal_H(M_{nP})$ is equal to
$\qfor_H(M)$.
\end{lem}

\textsc{Proof}. Theorem \ref{theo:Q-R-singular} - ``Quantization
commutes with reduction in the singular setting'' - tells us that
$\Qcal_H(M_{nP})= \sum_{\nu\in
\what{H}}\Qcal((M_{nP})_{\nu,H})V_\nu^H$ where $(M_{nP})_{\nu,H}$ is
the symplectic reduction
$$
(M_{nP}\times \overline{H\cdot\nu})/\!\!/_{0}H\cong(M\times
\Xcal_{nP}\times \overline{H\cdot\mu})/\!\!/_{(0,0)}H\times K.
$$

For $c>0$ we consider the compact subset of $\Kcal_c$ defined in
(\ref{eq:C-c}). Let $n_c\in\N$ such that $\Kcal_c\subset K\cdot (n
P^o)$ for any $n\geq n_c$. It implies that
$$
\Phi^{-1}_H\left(\xi\in\hgot^*, \|\xi\|\leq c\right)\subset
\Phi^{-1}_K(K\cdot (n P^o))
$$
for $n\geq n_c$. Since $M_{nP}$ "contains" as an the open subset
$\Phi^{-1}_K(K\cdot (n P^o))$, the arguments similar to those used
in the proof of Proposition \ref{prop:Q-M-P} show that
$\Qcal((M_{nP})_{\nu,H})=\Qcal(M_{\nu,H})$ for $\|\nu\|\leq c$ and
$n\geq n_c$. It means that
$$
\Qcal_H(M_{nP})= \sum_{\|\nu\|\leq c}\Qcal(M_{\nu,H})
    V_\nu^H + O_H(c)\quad \mathrm{when} \quad n\geq n_c.
$$
It follows that
$\lim\limits_{n\to\infty}\Qcal_H(M_{nP})=\sum_{\nu\in\what{H}}\Qcal(M_{\nu,H})
V_\nu^H=\qfor_H(M)$. $\Box$

\section{The case of an Hermitian space}

Let $(E,\h)$ be an Hermitian vector space of dimension $n$.

\subsection{The quantization of $E$}

Let
$\U:=\U(E)$ be the unitary group with Lie algebra $\ugot$.
We use the isomorphism $\epsilon:\ugot\to\ugot^*$ defined by
$\langle\epsilon(X),Y\rangle=-\Tr(XY)\in \R$. For $v,w\in E$, let
$v\otimes w^*: E\to E$ be the linear map $x\mapsto \h(x,w)v$.

Let $E_\R$ be the space $E$ view as a real vector space. Let $\Omega$ be
the imaginary part of $-\h$, and let $J$ the complex structure on
$E_\R$.  Then on $E_\R$, $\Omega$ is a (constant) symplectic
structure and $\Omega(-,J-)$ defines a scalar product. The action of
$\U$ on $(E_\R,\Omega)$ is Hamiltonian with moment map $\Phi:
E\to\ugot^*$ defined by
$\langle\Phi(v),X\rangle=\frac{1}{2}\Omega(Xv,v)$. Through $\epsilon$,
the moment map $\Phi$ is defined by
\begin{equation}\label{eq:Phi-E}
    \Phi(v)= \frac{1}{2i}v\otimes v^*.
\end{equation}

The pre-quantization data $(L,\langle -,-\rangle,\nabla)$ on the
Hamiltonian $\U$-manifold \break $(E_\R,\Omega,\Phi)$ is a trivial
line bundle $L$ with a trivial action of $\U$ equipped with the
Hermitian structure $\langle s,s'\rangle_v=
e^{\frac{-\h(v,v)}{2}}s\overline{s'}$ and the Hermitian connexion
$\nabla=d- i\theta$ where $\theta$ is the $1$-form on $E$ defined by
$\theta= \frac{1}{2}\Omega(v,dv)$.

\medskip

The traditional quantization of the Hamiltonian $\U$-manifold
$(E_\R,\Omega,\Phi)$, that we denote $\Qcal_U^{{\rm L}^2}(E)$, is
the {\em Bargman space} of entire holomorphic functions on $E$ which
are ${\rm L}^2$ integrable with respect to the Gaussian measure
$e^{\frac{-\h(v,v)}{2}}\Omega^n$. The representation $\Qcal_U^{{\rm L}^2}(E)$
of $\U$ is admissible. The irreducible representations of $\U$ that occur in
$\Qcal_U^{{\rm L}^2}(E)$ are the vector subspaces $S^j(E^*)$ formed by
the homogeneous polynomial on $E$ of degree $j\geq 0$.

On the other hand, the moment map $\Phi$ is proper (see
\ref{eq:Phi-E}). Hence we can consider the formal quantization
$\qfor_\U(E)\in\Rfor(\U)$ of the $\U$-action on $E$.

\medskip

\begin{lem}\label{lem:Q-formal-L2}
    The two quantizations of $(E,\Omega,\Phi)$,
$\Qcal_U^{{\rm L}^2}(E)$ and $\qfor_\U(E)$ coincide in $\Rfor(\U)$. In other
words, we have
\begin{equation}\label{eq:Q-formal-L2}
    \qfor_\U(E)= S^\bullet(E^*):=\sum_{j\geq 0} S^j(E^*) \quad \mathrm{in} \quad \Rfor(U).
\end{equation}
\end{lem}

\textsc{Proof}. Let $T\subset\U$ be a maximal torus with Lie algebra $\tgot\subset
\ugot$. There exists an
orthonormal basis $(e_k)_{k=1,\cdots,n}$ of $E$ and characters
$(\chi_k)_{k=1,\cdots,n}$ of $T$ such that $t\cdot e_k=\chi_k(t)e_k$
for all $k$. The family $(i e_k\otimes e_k^*)_{k=1,\cdots,n}$ is
then a basis of $\tgot$ such that $\frac{1}{i}d\chi_l(i e_k\otimes
e_k^*)=\delta_{l,k}$. The set $\what{U}\subset\tgot^*\subset\ugot^*$
of dominants weights
is composed, through $\epsilon$, by the elements
$$
\underline{\lambda}= i\sum_{k=1}^n \lambda_k e_k\otimes e_k^*,
$$
where $\lambda=(\lambda_1,\lambda_2,\cdots,\lambda_n)$ is a \emph{decreasing}
sequence of integer.

The formal quantization $\qfor_\U(E)\in\Rfor(\U)$ is defined
by
$$
\qfor_\U(E)=\sum_{\lambda_1\geq \cdots\geq \lambda_n} \Qcal(E_{\underline{\lambda}})
V_{\underline{\lambda}}
$$
where $E_{\underline{\lambda}}=\Phi^{-1}(U\cdot
\underline{\lambda})/U$ is the reduced space and
$V_{\underline{\lambda}}$ is the irreducible representation of $\U$
with highest weight $\underline{\lambda}$.

It is now easy to check that
$$
E_{\underline{\lambda}}=
\begin{cases}
    \{{\rm pt}\} & \text{if} \quad \lambda=(0,\cdots,0,-j)\  \mathrm{with}\  j\geq 0,\\
    \quad \emptyset        & \text{in \ the \ other \ cases} ,
\end{cases}
$$
and then
$$
\Qcal(E_{\underline{\lambda}})=
\begin{cases}
    1& \text{if} \quad \lambda=(0,\cdots,0,-j)\  \mathrm{with}\  j\geq 0,\\
    0     & \text{in \ the \ other \ cases} .
\end{cases}
$$
Finally (\ref{eq:Q-formal-L2}) follows from the fact that
$V_{\underline{(0,\cdots,0,-j)}}= S^j(E^*)$. $\Box$

\subsection{The quantization of $E$ restricted to a subgroup of $\U$}

Let
$K\subset\U$ be a connected Lie subgroup with Lie algebra $\kgot^*$. Let
$K_\C\subset {\rm GL}(E)$ be its complexification. The moment map
relative to the $K$-action on $(E_\R,\Omega)$ is the map
$$
\Phi_K: E\to \kgot^*
$$
equal to the composition of $\Phi$ with the projection $\ugot^*\to\kgot^*$.

\begin{lem}\label{lem:closed-orbit-E}
    The following conditions are equivalent :
    \begin{itemize}
      \item[(a)] the map $\Phi_K$ is proper,
      \item[(b)] $\Phi_K^{-1}(0)=\{0\}$,
      \item[(c)] $\{0\}$ is the only closed $K_\C$-orbit in $E$,
      \item[(d)] for every $v\in E$ we have $0\in \overline{K_\C\cdot v}$,
      \item[(e)] $S^\bullet(E^*)$ is an admissible representation of $K$,
      \item[(f)] the $K$-invariant polynomials on $E$ are the constant polynomials.
    \end{itemize}
\end{lem}

\textsc{Proof}. The equivalence $(a)\Longleftrightarrow (b)$ is due to the fact that
$\Phi_K$ is quadratic.

Let $\Ocal$ be a $K_\C$-orbit in $E$. Classical results of Geometric
Invariant Theory \cite{GIT,Kempf-Ness} assert that
$\overline{\Ocal}\cap\Phi_K^{-1}(0)\neq \emptyset$ and that $\Ocal$
is closed if and only if $\Ocal\cap\Phi_K^{-1}(0)\neq \emptyset$.
Hence $(b)\Longleftrightarrow (c)\Longleftrightarrow (d)$.

After Lemma \ref{lem:Q-formal-L2} we know that $\qfor_\U(E)=
S^\bullet(E^*)$. Since $\qfor_\U(E)$ is $K$-admissible when $\Phi_K$
is proper (see Section \ref{sec:proof}), we have $(a)\Longrightarrow
(e)$.

For every $\mu\in \what{K}$, the $\mu$-isotopic component
$[S^\bullet(E^*)]_\mu$ is a module over
$[S^\bullet(E^*)]_0=[S^\bullet(E^*)]^K$. Hence $\dim
[S^\bullet(E^*)]_\mu < \infty$ implies that $[S^\bullet(E^*)]^K=\C$.
We have $(e)\Longrightarrow (f)$.

Finally $(f)\Longrightarrow (d)$ follows from the following
fundamental fact. For any $v,w\in E$ we have $\overline{K_\C\cdot
v}\cap \overline{K_\C\cdot w}\neq \emptyset$ if and only if
$P(v)=P(w)$ for all $P\in[S^\bullet(E^*)]^K$.  $\Box$

\bigskip

Theorem \ref{theo:intro} implies the following

\begin{prop}\label{prop-E-K-mu}
Let $K\subset \U(E)$ be a closed connected subgroup such that
$S^\bullet(E^*)$ is an admissible representation of $K$.
For every $\mu\in \what{K}$, we have
$$
\dim\left([S^\bullet(E^*)]_\mu\right)=\Qcal(E_{\mu,K})
$$
where $[S^\bullet(E^*)]_\mu$ is the $\mu$-isotopic component of $S^\bullet(E^*)$
and $E_{\mu,K}$ is the reduced space
$\Phi_K^{-1}(K\cdot \mu)/K$.
\end{prop}

\medskip

In the following examples the condition $\Phi_K^{-1}(0)=\{0\}$ is
easy to check.

\begin{itemize}
      \item[1)] the subgroup $K\subset \U(E)$ contains the center of $\U(E)$,
      \item[2)] $E=\wedge^2 \C^n$ or $E=S^2(\C^n)$ and $K=\U(n)\subset \U(E)$,
      \item[3)] $E={\rm M}_{n,k}$ is the vector space of $n\times k$-matrices
      and $K=\U(n)\times \U(k)\subset \U(E)$.
\end{itemize}



\end{document}